\documentclass[12pt]{article}
\usepackage{amsfonts}
\usepackage{amssymb}
\usepackage{polski}

\input{tcilatex}
\begin{document}

\begin{center}
\textbf{Beurling moving averages and approximate homomorphisms}

\textbf{by \\[0pt]
N. H. Bingham and A. J. Ostaszewski}\\[0pt]

\bigskip

\bigskip
\end{center}

\noindent \textbf{Abstract. }The theory of regular variation, in its
Karamata and Bojani\'{c}-Karamata/de Haan forms, is long established and
makes essential use of homomorphisms. Both forms are subsumed within the
recent theory of Beurling regular variation, developed further here,
especially certain moving averages occurring there. Extensive use of group
structures leads to an algebraicization not previously encountered here, and
to the approximate homomorphisms of the title. Dichotomy results are
obtained: things are either very nice or very nasty. Quantifier weakening is
extended, and the degradation resulting from working with limsup and liminf,
rather than assuming limits exist, is studied.

\bigskip


\bigskip


\begin{center}
\textsc{CONTENTS}
\end{center}

{\small 1. Introduction}

{\small 2. From Beurling to Karamata}

{\small 3. Popa (circle) groups}

{\small 4. Extensions to Beurling's Tauberian Theorem}

{\small 5. Uniformity, semicontinuity}

{\small 6. Dichotomy}

{\small 7. Quantifier weakening}

{\small 8. Representation}

{\small 9. Divided difference and double sweep}

{\small 10. Uniform Boundedness Theorem}

{\small 11. Character degradation from limsup}

{\small References}

{\small Appendix}\textsl{\newline
}

\section{Introduction}

This work is a sequel to our recent papers [BinO10,11,12] together with the
related papers [Ost2,3,4] by the second author and one [Bin] by the first,
re-examined in the light of two much earlier works [BinG2,3] by the first
author and Goldie. Our title Beurling moving averages addresses both the
Beurling slow and regular variation in [BinO10] (to which we refer for
background), and [BinG2,3], the motivation for which is strong laws of large
numbers in probability theory.

Beurling regular variation is closely linked with Karamata regular variation
(the standard work on which is [BinGT], BGT\ below, to which we refer for
background). In [BinO10], it emerged that Beurling regular variation in fact
subsumes the traditional (and very widely used) Karamata regular variation,
together with its Bojani\'{c}-Karamata/de Haan relative -- BGT Ch. 1-3;
[BojK], [dH]. Whereas the traditional approach is to develop the measure and
Baire-property (`Baire' below) cases in parallel, measure being regarded as
primary, it is now clear both that one can subsume both cases together and
that it is in fact the Baire case that is primary. This is the theory of
topological regular variation, for which see [BinO1,2,4,5], [Ost1] -- this
informs our approach in \S 10.

It is convenient to work both multiplicatively in $\mathbb{R}_{+}:=(0,\infty
)$ and additively in $\mathbb{R}.$ A self-map $f$ of $\mathbb{R}_{+}$ or $h$
of $\mathbb{R}$ is \textit{Beurling }$\varphi $\textit{-slowly varying} if,
according to context,
\begin{equation}
f(x+t\varphi (x))/f(x)\rightarrow 1,\text{ or }h(x+u\varphi
(x))-h(x)\rightarrow 0,  \tag{$BSV/BSV_{+}$}
\end{equation}%
as $x\rightarrow \infty ,$ where $\varphi $ is a self-map of $\mathbb{R}_{+}$
and is \textit{self-neglecting} ($\varphi \in SN$), so that
\begin{equation}
\varphi (x+t\varphi (x))/\varphi (x)\rightarrow 1\text{ locally uniformly in
}t\text{ for all }t\in \mathbb{R}_{+},  \tag{$SN$}
\end{equation}%
and $\varphi (x)=o(x).$ This traditional restriction may be usefully relaxed
in two ways, as in [Ost3]: firstly, in imposing the weaker order condition $%
\varphi (x)=O(x),$ and secondly by replacing the limit $1$ by a general
limit function $\eta >0,$ so that for $\mathbb{A}=[0,\infty )$
\begin{equation}
\eta _{x}^{\varphi }(t):=\varphi (x+t\varphi (x))/\varphi (x)\rightarrow
\eta (t)>0\text{ locally uniformly in }t\text{ for }t\in \mathbb{A}.
\tag{$SE_{\QTR{Bbb}{A}}$}
\end{equation}%
Following [Ost3], such a $\varphi $ will be called \textit{self-equivarying}%
, $\varphi \in SE,$ and the limit function\footnote{%
Note the changes here: positivity has been incorporated into the definition
(for simplicity), $\eta ^{\varphi }$ replaces the original notation $\lambda
_{\varphi }$ for this context, both to free up the use of $\lambda $ for
other conventional uses, and to reflect the connection to the function $%
H_{\rho }$ below (as $H$ denotes the Greek capital `eta'). Finally, $t=0$ is
included under $(SE_{\mathbb{A}})$ above, being a consequence of the
assertion for $t>0$ -- see \S 5 Lemma 1, Theorem 3.} $\eta =\eta ^{\varphi }$
necessarily satisfies the \textit{Beurling functional equation}%
\begin{equation}
\eta (u+v\eta (u))=\eta (u)\eta (v)\text{ for }u,v\in \mathbb{R}_{+}
\tag{$BFE$}
\end{equation}%
(this is a special case of the \textit{Go\l \k{a}b-Schinzel equation} $(GS),$
here conditioned by its relation to $(SE_{\mathbb{A}})$ -- see also e.g.
[Brz1], [BrzM], or [BinO11]). As $\eta \geq 0,$ imposing the natural
condition $\eta >0$ (on $\mathbb{R}_{+}$) above implies that it is
continuous and of the form%
\[
\eta (t)\equiv \eta _{\rho }(t):=1+\rho t,\quad (t\geq 0)\qquad \text{for
some }\rho \geq 0
\]%
(see [BinO11]). Then we call $\eta $ a \textit{Beurling function: }$\eta \in
GS,$ with $\rho $ the $\eta $\textit{-index} (of $\varphi $ when $\eta =\eta
^{\varphi },$ and then we write $\rho =\rho _{\varphi }$); as in BGT 2.11,
we extend in \S 5 the domain (and local uniformity in $(SE_{\mathbb{A}})$)
to $\mathbb{A}=(\rho ^{\ast },\infty ),$ where $\rho ^{\ast }:=-\rho ^{-1}$;
in \S 3 we call $\rho ^{\ast }$ the \textit{Popa origin.} The case $\rho =0$
recovers $SN$. For $\varphi \in SE,$ a self-map $f$ of $\mathbb{R}_{+}$ or $%
h $ of $\mathbb{R}$ is \textit{Beurling }$\varphi $\textit{-regularly varying%
} if, according to context, the limits below exist:
\begin{equation}
f(x+t\varphi (x))/f(x)\rightarrow g(t),\text{ or }h(x+u\varphi
(x))-h(x)\rightarrow k(u).  \tag{$BRV/BRV_{+}$}
\end{equation}%
For $\varphi \in SN$ and $f$ Baire/measurable, the limit $g(t)$ is
necessarily an exponential function $e^{\gamma t}$ (provided $g>0$ on a
non-negligible set), equivalently $k$ is linear: $k(u)\equiv \gamma u$,
convergence is locally uniform, and there is a representation for the
possible $f$ (see [BinO10]), involving the $\varphi $\textit{-index of
Beurling variation,} or \textit{Beurling }$\varphi $\textit{-index} for
short, $\gamma $. For $\varphi \in SE$ with $\eta $-index $\rho >0$, the
situation is altered from $g(t)=e^{\gamma t}$ so that (see [Ost3, Th. 1$%
^{\prime }$])%
\begin{equation}
g(t)=(1+\rho t)^{\gamma },\text{ or }k(t)=\gamma \log (1+\rho t)\qquad
(t>\rho ^{\ast }).  \tag{$\rho $ -$BR_{\gamma }$}
\end{equation}%
\indent We are led to the question of existence and additivity properties of
the limit functions below:%
\[
K_{F}(t):=\lim_{x\rightarrow \infty }\left. \Delta _{t}^{\varphi
}F(x)\right/ \varphi (x),\text{ }K_{F}^{\ast }(t):=\limsup_{x\rightarrow
\infty }\left. \Delta _{t}^{\varphi }F(x)\right/ \varphi (x),
\]%
with $\Delta _{t}^{\varphi }$ the difference operator%
\[
\Delta _{t}^{\varphi }F(x):=F(x+t\varphi (x))-F(x),
\]%
and local uniform convergence assumed (unless otherwise stated). For $%
\varphi (x)\equiv 1$ this reduces to the usual difference operator $\Delta
_{t}$. Motivated by classical analysis, we introduce a more general
auxiliary function $\psi (x)$ in the denominator:%
\[
K_{F}(t):=\lim \left. \Delta _{t}^{\varphi }F(x)\right/ \psi (x),\text{ }%
K_{F}^{\ast }(t):=\lim \sup \left. \Delta _{t}^{\varphi }F(x)\right/ \psi
(x).
\]%
If $K_{F}$ is defined at $u$ and $v$, then (cf. \S 8 Lemma 3)
\[
K_{F}(v+uh(v))=K_{F}(v)+K_{F}(u)g(v),
\]%
provided%
\[
h(v):=\lim \left. \varphi (x+v\varphi (x))\right/ \varphi (x)\text{ and }%
g(u):=\lim \left. \psi (x+u\varphi (x))\right/ \psi (x)
\]%
exist (and convergence to $K_{F}$ is locally uniform), which will be the
case when $\varphi \in SE$ (so that $h=\eta _{\rho })$ and $\psi $ is $%
\varphi $-regularly varying (so that either $\rho =0$ and $g=e^{\gamma \cdot
},$ or $\rho >0$ and $g\equiv (1+\rho \cdot )^{\gamma },$ by $(\rho $ -$%
BR_{\gamma })$ above). The related functional equation -- the extended
Goldie-Beurling (Pexiderized\footnote{%
After Pexider's equation: $f(xy)=g(x)+h(y)$ in three unknown functions and
its generalizations -- cf. [Kuc, 13.3], [Brz1, 2]. See also [Ste] for the
more general Levi-Civita functional equations.}) equation,%
\begin{equation}
K(v+uh(v))=K(v)+\kappa (u)g(v),  \tag{$GBE$-$P$}
\end{equation}%
for $h,\kappa $ positive -- is studied in [BinO11, Th. 9 and 10]; special
cases appear below in \S 2 Cor. 2, \S 8 Lemma 3, \S 9 Prop. 10. Its
solutions $K$, necessarily continuous, are there characterized (subject to $%
K(0)=0$) as%
\[
K(x)\equiv c\cdot \tau _{f}(x)\text{ with }f:=h/g\text{ and }\tau
_{f}(x):=\int_{0}^{x}\mathrm{d}w/f(w)\qquad (x\geq 0),
\]%
an `occupation time measure' (of the interval $[0,x];$ \S 2) and $c\in
\mathbb{R}$; the `relative flow rate' $f$ satisfies the \textit{%
Cauchy-Beurling exponential equation}:%
\begin{equation}
f(v\circ _{h}u)=f(u)f(v),  \tag{$CBE$}
\end{equation}%
cf. [Chu], [Ost4]. Here $\circ _{h}$ denotes Popa's binary operation ([Pop],
cf. [Jav], \S 3 below)
\[
v\circ _{h}u:=v+uh(v),
\]%
so that $h=\eta _{\rho }$ itself also satisfies $(CBE)$; this confers a
group structure, turning certain subsets of $\mathbb{R}$ into groups, called
\textit{Popa (circle) groups} in \S 3; furthermore, necessarily $\kappa =K$.
Solving $(GBE$-$P)$ may be expressed as an equivalent Popa \textit{%
homomorphism problem} of finding $k,h\in GS$ satisfying
\begin{equation}
K(v\circ _{h}u)=K(v)\circ _{k}K(u)  \tag{$GBE$}
\end{equation}%
(cf. [Brz2], [Mur], [Ost4]), where%
\[
k(u)=g(K^{-1}(u)).
\]%
This observation is new even for the classical context $h\equiv 1;$ here $%
f=e^{-\gamma t},$ so
\[
\tau _{f}(x)\equiv H_{\gamma }(x):=(e^{\gamma x}-1)/\gamma \text{ with }%
H_{0}(x)\equiv x.
\]%
For $\eta \equiv \eta _{\rho }$ with $\rho >0,$ $g\equiv (1+\rho \cdot
)^{\gamma },$ by $(\rho $ -$BR_{\gamma })$ above, $f(x)=(1+\rho x)^{1-\gamma
},$ so for $x>\rho ^{\ast },$%
\[
K\equiv c\cdot \tau _{f}=c\cdot K_{\rho \gamma },\text{ where }K_{\rho
\gamma }(x):=\int_{0}^{x}(1+\rho w)^{\gamma -1}\mathrm{d}w=\left( (1+\rho
x)^{\gamma }-1\right) /\rho \gamma
\]%
(linear for $\gamma =1$). The `slow case' $\gamma =0$ may also be handled via%
\[
\lim\nolimits_{\gamma \rightarrow 0}K_{\rho \gamma }(x)=\log (1+\rho x)/\rho
\qquad (x>\rho ^{\ast }).
\]%
\newline
\textit{\noindent }When $\varphi (x)\equiv 1,$ the moving averages $\left.
\Delta _{t}^{\varphi }F(x)\right/ \psi (x)$ reduce to classical Bojani\'{c}%
-Karamata/de Haan limits (BGT\ Ch. 3), for which the auxiliary $\psi (x)$ is
necessarily Karamata regularly varying, so just as before (trivially, since $%
\varphi \in SE)$ has exponential limit function, $g\equiv e^{\gamma \cdot }$
say, and then $(GBE$-$P)$ simplifies to the original \textit{Goldie
functional equation} (see e.g. [BinO11,12], [Ost4]):%
\begin{equation}
K(u+v)=e^{\gamma u}K(v)+K(u),  \tag{$GFE$}
\end{equation}%
with solution $K(u)\equiv c\cdot H_{\gamma }(u),$ as before. The latter
function plays a crucial role in the Bojani\'{c}-Karamata/de Haan theory of
regular variation. Here, and in the general case, if $\Delta _{t}^{\varphi
}F/\psi $ has a limiting moving average $K_{F}$, then for some $c_{F}\in
\mathbb{R}$, as above (cf. [BinO11, Th. 3, 9, 10]),
\[
K_{F}(u)=c_{F}\cdot H_{\gamma }(u),
\]%
with $c_{F}$ the $\psi $\textit{-index} of $F$ (for $\psi $ which is $%
\varphi $-regularly varying), while $\psi $ has \textit{Beurling }$\varphi $%
\textit{-index} $\gamma .$

In the classical context, with \textrm{limsup} in place of limit one works
also with $K_{F}^{\ast },$ abbreviated to $K^{\ast }$ (and similarly $%
K_{\ast }$ with \textrm{liminf}). Here the equations $(GFE)\ $give way to
functional inequalities, such as the \textit{Goldie functional inequality}
\begin{equation}
K^{\ast }(u+v)\leq e^{\gamma u}K^{\ast }(v)+K^{\ast }(u)  \tag{$GFI$}
\end{equation}%
(BGT (3.2.5)), which we summarize by saying that $K^{\ast }$ is \textit{%
exp-subadditive. }Equivalently, this may be re-expressed symmetrically here
as \textit{group sub-additivity}:%
\[
K^{\ast }(x+y)\leq K^{\ast }(x)\circ _{k}K^{\ast }(y)
\]%
with $k$ as above, and in the more general Beurling case correspondingly to $%
(GBE)$ as%
\[
K^{\ast }(x\circ _{h}y)\leq K^{\ast }(x)\circ _{k}K^{\ast }(y).
\]

For $\psi $ regularly varying, the set%
\[
\mathbb{A}:=\{t:\lim \Delta _{t}F(x)/\psi (x)\text{ exists and is finite}\},
\]%
for which see e.g. BGT Th. 3.2.5 (proof) and \S \S 5,6 below, constitutes
the domain of the function%
\begin{equation}
K_{F}(a):=\lim_{x\rightarrow \infty }\Delta _{a}F(x)/\psi (x)\qquad (a\in
\mathbb{A});  \tag{ker}
\end{equation}%
hence we refer to $K_{F}$ here and above as the \textit{regular kernel }of $%
F $ -- the homomorphism approximating $F$ of our title. In [BinO11] (and in
[BinO12] for the case $\rho =0),$ we study conditions on $K^{\ast }$
implying that $K_{F}$ exists, i.e. that the inequality becomes an equation,
by imposing `Heiberg-Seneta' side-conditions (see \S 7 Prop. 9), and density
of $\mathbb{A}$ -- again cf. BGT Ch. 3, especially the crucial Theorem
4.2.5. Below these findings are extended to the Beurling context.

In view of the algebraic treatment to follow in \S 3 on Popa groups, one may
regard the terms \textit{additive }and \textit{homomorphic }as synonymous
for our purposes here.

\section{From Beurling to Karamata}

The function $H_{\rho }$ (of \S 1) satisfies%
\[
\left. dH_{\rho }\right/ dx=e^{\rho x}=1+\rho H_{\rho }(x)=\eta _{\rho
}(H_{\rho }(x)),
\]%
and solves the Goldie equation $(GFE),$ in which the auxiliary function $g,$
which is necessarily exponential for $K$ Baire/measurable, takes the form $%
g(x)=e^{\rho x}$ -- again see [BinO11, Th. 1]. Regarding $\varphi ,\eta \in
SE$ as generating (velocity) flows as in [BinO10], their occupation `times'
(on $[0,x])$ are (cf. [Bec, p.153]):%
\[
\tau _{\varphi }(x):=\int_{0}^{x}\mathrm{d}w/\varphi (w)\text{ and }\tau
_{\eta }(x):=\int_{0}^{x}\mathrm{d}w/\eta (w),
\]%
both strictly increasing. (For present needs this notation is more
symmetrical than that of [BinG1] with $\Phi $ for $\tau _{\varphi },$ and of
BGT 2.12.29, which we mention for purposes of comparison.) For $\rho >0$ and
$\eta =\eta _{\rho }\in SE$
\[
\tau _{\eta }(x):=\int_{0}^{x}\frac{\mathrm{d}w}{1+\rho w}=\frac{1}{\rho }%
\log (1+\rho x),
\]%
so
\[
\tau _{\eta }^{-1}(t)=H_{\rho }(t)=(e^{\rho t}-1)/\rho .
\]%
In particular, the trajectory $w(t):=\tau _{\eta }^{-1}(t)$ satisfies the
equation
\[
\mathrm{d}w(t)/\mathrm{d}t=e^{\rho t}=1+\rho w(t)=\eta (w(t))\text{ with }%
w(0)=0.
\]%
Necessarily, working with the (inverse) re-parametrization $\mathrm{d}t(w)/%
\mathrm{d}w=e^{-\rho t}=\psi (t)\in SE$ gives $\tau _{\psi }(x)=H_{\rho
}(x), $ again an occupation time measure.

We now generalize a theorem of Bingham and Goldie [BinG2, Th. 2]. This
recovers their theorem when $\rho _{\eta }=0$ and $\varphi (x)=o(x),$ as
then $\varphi \in SN.$ The result may be interpreted as a local `chain
rule', for $V(s)=U(s(t)),$ where the trajectory $s(t):=\tau _{\varphi
}^{-1}(t)$ satisfies $\mathrm{d}s(t)/\mathrm{d}t=\varphi (s(t))=\varphi
(\tau _{\varphi }^{-1}(t))=g(t)$ (with $\varphi \in SE,$ a `self-equivarying
flow').

\bigskip

\noindent \textbf{Theorem 1 (Time-change Equivalence Theorem for Moving
Averages).} \textit{For positive }$\varphi \in SE$\textit{\ with }$1/\varphi
$\textit{\ locally integrable, }$U$\textit{\ satisfies}%
\begin{equation}
\frac{U(x+t\varphi (x))-U(x)}{\varphi (x)}\rightarrow c_{U}t\text{ as }%
x\rightarrow \infty ,\text{ for all }t\geq 0  \tag{$BMA_{\varphi }$}
\end{equation}%
\textit{iff its \textbf{time-changed version} }$V:=U\circ \tau _{\varphi
}^{-1}$\textit{\ satisfies, for }$g(y):=\varphi (\tau _{\varphi }^{-1}(y)),$%
\textit{\ }%
\begin{equation}
\frac{V(y+s)-V(y)}{g(y)}\rightarrow c_{U}H_{\rho }(s)\text{ as }y\rightarrow
\infty ,\text{ for all }s\geq 0,  \tag{$KMA_{g}$}
\end{equation}%
\textit{where }$\rho =\rho _{\varphi }$\textit{\ is the }$\eta $\textit{%
-index of }$\varphi $\textit{.}

\bigskip

This is proved exactly as in [BinG2, Th. 2], using the following.

\bigskip

\noindent \textbf{Proposition 1. }\textit{For }$\varphi \in SE$\textit{\ and
}$\eta =\eta ^{\varphi }$, \textit{locally uniformly in }$s$%
\[
\lim [\tau _{\varphi }(x+s\varphi (x))-\tau _{\varphi }(x)]=\tau _{\eta
}(s).
\]%
\textit{In particular, this is so for }$\varphi \in SN$, \textit{where }$%
\tau _{\eta }(s)\equiv s.$

\bigskip

\noindent \textit{Proof.}\textbf{\ }Let\textbf{\ }$\rho $ be the $\eta $%
-index. Fix $s>0,$ then uniformly in $t\in \lbrack 0,s]$%
\[
\varepsilon (x,t):=\varphi (x)/\varphi (x+t\varphi (x))-1/\eta
(t)\rightarrow 0,\text{ so }e(x,s):=\int_{0}^{s}\varepsilon (x,t)\mathrm{d}%
t\rightarrow 0.
\]%
Then, as in [BinG2, Th. 2], using the substitution $w=x+t\varphi (x)$%
\begin{eqnarray*}
\tau _{\varphi }(x+s\varphi (x))-\tau _{\varphi }(x) &=&\int_{x}^{x+s\varphi
(x)}\mathrm{d}w/\varphi (w)=\int_{0}^{s}\frac{\varphi (x)\mathrm{d}t}{%
\varphi (x+t\varphi (x))} \\
&=&\int_{0}^{s}\left( \frac{1}{\eta (t)}+\varepsilon (x,t)\right) \mathrm{d}%
t=\tau _{\eta }(s)+e(x,u).\text{ }
\end{eqnarray*}%
If $\varphi \in SN,$ then $\tau _{\eta }(s)\equiv s,$ as $\eta \equiv 1.$ $%
\square $

\bigskip

Our first corollary characterizes $SE$ in terms of a \textit{multiplicative}
Karamata index via its time-changed version $g$; this is a \textit{%
consistency} result in view of the characterization from [Ost3] of $\varphi
\in SE$ as the product $\eta ^{\varphi }\psi $ with $\psi $ in $SN.$ The
latter identifies $\varphi $ itself as having \textit{additive} Karamata
index $\rho _{\varphi }.$

\bigskip

\noindent \textbf{Corollary 1. }$\varphi \in SE$\textit{\ iff }$g=\varphi
\circ \tau _{\varphi }^{-1}$ \textit{is regularly varying in the
additive-argument sense with multiplicative Karamata index }$\rho _{\varphi
}.$ \textit{In particular, }$\varphi \in SN$\textit{\ iff }$g=\varphi \circ
\tau _{\varphi }^{-1}$ \textit{is regularly varying with multiplicative
Karamata index }$\rho _{\varphi }=0.$

\bigskip

\noindent \textit{Proof.}\textbf{\ }Put $\rho =\rho _{\varphi }.$ Since%
\[
\left. (\varphi (x+t\varphi (x))-\varphi (x))\right/ \varphi (x)=\left.
\varphi (x+t\varphi (x))\right/ \varphi (x)-1\rightarrow \rho t,
\]%
we may apply Th. 0 to $U=\varphi $ so that $V:=\varphi \circ \tau _{\varphi
}^{-1}=g;$ then by $(KMA_{g})$%
\[
\left. g(y+s)\right/ g(y)-1=\left. (g(y+s)-g(y))\right/ g(y)\rightarrow
(e^{\rho _{\varphi }x}-1):\quad \left. g(y+s)\right/ g(y)\rightarrow e^{\rho
x},
\]%
and conversely. $\square $

\bigskip

If $K_{V}(s)$ -- defined by (ker) above (with $g$ for $\psi $) -- exists for
all $s$, as in $(KMA_{g})$, then as we now show $K_{V}$ satisfies a Goldie
equation, from which the form of $K_{V}$ can be read off, as in the
Equivalence Theorem, Theorem 1.

\bigskip\

\noindent \textbf{Corollary 2. }\textit{For }$\varphi \in SE,$\textit{\ so
that }$g=\varphi \circ \tau _{\varphi }^{-1}$ \textit{is regularly varying
with multiplicative Karamata index }$\rho =\rho _{\varphi }:$

\noindent \textit{if }$KMA_{g}$ \textit{\ -- equivalently }$BMA_{\varphi }$
\textit{-- holds, then for }$K_{V}(u),$ \textit{as above,}%
\[
K_{V}(s+t)=K_{V}(s)e^{\rho t}+K_{V}(t),
\]%
\textit{and so for some }$c$%
\[
K_{V}(s)=cH_{\rho }(s).
\]

\noindent \textit{Proof.} The Goldie equation follows from Corollary 1, since%
\[
\frac{V(y+s+t)-V(y)}{g(y)}=\frac{V(y+s+t)-V(y+t)}{g(y+t)}\frac{g(y+t)}{g(y)}+%
\frac{V(y+t)-V(y)}{g(y)}.
\]%
Now apply Theorem 2 of [BinO11] to deduce the form of $K_{V}.$ $\square $

\section{Popa (circle) groups}

Recall from Popa [Pop], for $h:\mathbb{R}\rightarrow \mathbb{R}$, the
\textit{Popa operation} $\circ _{h}$ and its \textit{Popa domain} $\mathbb{G}%
_{h}$ (our terminology) defined by:%
\[
a\circ _{h}b:=a+bh(a),\qquad \text{ }\mathbb{G}_{h}:=\{g:h(g)\neq 0\}.
\]%
The special case (but nevertheless typical -- see below) of $h(t)=\eta
_{1}(t)\equiv 1+t$ yields the \textit{circle product} in a ring, $a\circ
b:=a+b+ab$ -- see [Ost4] for background. We recall also, from Javor [Jav]
(in the broader context of $h:\mathbb{E}\rightarrow \mathbb{F}$, with $%
\mathbb{E}$ a vector space over a commutative field $\mathbb{F}$), that $%
\circ _{h}$ is associative iff $h$ satisfies the Go\l \k{a}b-Schinzel
equation, briefly $h\in GS$ (cf. \S 1 -- a temporary ambiguity resolved
below):%
\begin{equation}
h(x+yh(x))=h(x)h(y)\qquad (x,y\in \mathbb{G}_{h}).  \tag{$GS$}
\end{equation}%
Their role below is fundamental; first, $GS\subseteq SE,$ and for $\varphi
\in SE$ the Popa operation $x\circ _{\varphi }t=x+t\varphi (x)$ compactly
expresses the Beurling transformation $t\rightarrow x+t\varphi (x).$ More is
true: taking one step further beyond $GS$ to $SE$ is an operation localized
to $x:$
\[
s\circ _{\varphi x}t:=s+t\eta _{x}^{\varphi }(s),\text{ where }\eta
_{x}^{\varphi }(s),\text{ or just }\eta _{x}(s):=\varphi (x+s\varphi
(x))/\varphi (x)
\]%
as in \S 1 (we use $\eta _{x}^{\varphi }$ or $\eta _{x}$ depending on
emphasis or context). The notation above neatly summarizes two frequently
used facts in (Karamata/Beurling) regular variation:%
\[
x\circ _{\varphi }(b\circ _{\varphi x}a)=y\circ _{\varphi }a,\text{ for }%
y=x+b\varphi (x)
\]%
(proved in Prop. 2(ii) below), and as $x\rightarrow \infty ,$ locally
uniformly in $s,t:$%
\[
s\circ _{\varphi x}t\rightarrow s\circ _{\eta }t,\text{ for }\eta
(s):=\lim\nolimits_{x}\eta _{x}^{\varphi }(s)\in GS.
\]%
So here we return to $GS.$

The appearance of a group structure `in the limit' is not accidental -- see
[Ost4] for background. The fact that, for $\eta $ as here, $\eta \in GS$ is
proved in [Ost3] -- see \S 1; solutions of $(GS)$ that are \textit{positive}
on $\mathbb{R}_{+}:=(0,\infty )$ are key here, being of the form $\eta
_{\rho }(x):=1+\rho x$ with $\rho \geq 0.$ The case $\rho =0$ corresponds to
the classical Karamata setting, and $\rho >0$ to the recently established,
general, theory of Beurling regular variation [BinO10]. For the
corresponding \textit{Popa groups} write $\circ _{\rho }$ (when $h=\eta
_{\rho }),$ or even $\circ ,$ omitting subscripts both on $\circ $ and on $%
\eta ,$ if context permits. To prevent confusion, $u_{\circ }^{-1}$ denotes
the relevant group inverse. Furthermore, we employ the notation:%
\[
\rho ^{\ast }:=-\rho ^{-1},\quad \mathbb{G}_{\ast }^{\rho }:=\mathbb{R}%
\backslash \{\rho ^{\ast }\},\quad \mathbb{G}_{+}^{\rho }:=(\rho ^{\ast
},\infty ),\quad \mathbb{G}_{-}^{\rho }:=(-\infty ,\rho ^{\ast }),\qquad
(\rho \neq 0),
\]%
\[
\begin{array}{cc}
\mathbb{G}_{\ast }^{\infty }:=\mathbb{R}^{\ast }=\mathbb{R}\backslash \{0\},
& \mathbb{G}_{\ast }^{0}:=\mathbb{R}, \\
\eta _{\ast }^{\rho }(x):=\eta _{\rho }(x)\text{ }(\rho \neq 0), & \eta
_{\ast }^{0}(x):=e^{x}.%
\end{array}%
\]%
We call $\rho ^{\ast }$ the \textit{Popa origin} (of $\mathbb{G}_{\ast
}^{\rho }$), interpreting it when $\rho =0$ as $-\infty .$ Its critical role
for Beurling regular variation emerges in \S 5 Lemma 1.

We collect relevant facts in the following, slightly extending work of Popa
[Pop, Prop. 2] and Javor [Jav, Lemma 1.2].

\bigskip

\noindent \textbf{Theorem PJ.}\textit{\ For } $\varphi $ \textit{satisfying }%
$(GS)$\textit{\ above, not the zero map, }$(\mathbb{G}_{\varphi },\circ
_{\varphi })$\textit{\ is a group. If }$\varphi $\textit{\ is injective on }$%
\mathbb{G}_{\varphi },$ \textit{then }$\circ _{\varphi }$\textit{\ is
commutative, and }$\varphi $ \textit{maps homomorphically into }$(\mathbb{R}%
^{\ast },\cdot )$:%
\[
\varphi (x\circ _{\varphi }y)=\varphi (x)\varphi (y).
\]%
\textit{In particular, }$\mathbb{G=G}^{\rho }:=(\mathbb{G}_{\ast }^{\rho
},\circ _{\rho })$\textit{\ is an abelian group with }$1_{\mathbb{G}}=0$
\textit{and inverse}%
\[
u_{\circ }^{-1}=-u/\eta _{\rho }(u).
\]%
$\mathbb{G}^{0}:=(\mathbb{R},\circ )$ \textit{is }$(\mathbb{R},+)$\textit{\
for }$\rho =0,$\textit{\ so that }$\mathbb{G}^{\rho }$\textit{\ is
isomorphic under }$\eta _{\ast }^{\rho }$ \textit{to} $(\mathbb{R}^{\ast
},\cdot )$ \textit{for }$\rho \geq 0$. \textit{Furthermore,}

\noindent (i) \textit{inversion carries }$\mathbb{G}_{+}^{\rho }$\textit{\
into itself: }$(\mathbb{G}_{+}^{\rho })_{\circ }^{-1}=\mathbb{G}_{+}^{\rho }$
\textit{and }$\eta _{\rho }^{\ast }$ \textit{carries }$\mathbb{G}_{+}^{\rho
} $\textit{\ onto }$\mathbb{R}_{+}$\textit{;}

\noindent (ii) \textit{for }$\rho >0$ \textit{the reflection }$\pi =\pi
_{\rho }:u\mapsto -u+2\rho ^{\ast }$ \textit{carries each of }$\mathbb{G}%
_{+}^{\rho }$ \textit{and }$\mathbb{G}_{-}^{\rho }$ \textit{%
skew-isomorphically onto the other in the sense that}%
\[
\pi ^{-1}(\pi (s)\circ _{\rho }\pi (t))=\pi (s\circ _{\rho }t),
\]%
\textit{\ and}%
\[
|\eta _{\ast }^{\rho }(\pi (t))|=\eta _{\ast }^{\rho }(t)\qquad (t\in
\mathbb{G}_{+}^{\rho });\qquad \eta _{\ast }^{\rho }(\rho ^{\ast })=0.
\]

\bigskip

\noindent \textit{Proof.} In general, if $\varphi $ is injective on $\mathbb{%
G}_{\varphi },$ then $\circ _{\varphi }$ is commutative, as $(GS)$ is
symmetric on the right-hand side. Commutativity of $\circ _{\rho }$ follows
directly from $v+u(1+\rho v)=u+v(1+\rho u).$ As $u\circ 0=u$ and $0\circ
v=v, $ the neutral element is $1_{\mathbb{G}}=0$; the inverse is%
\[
v_{\circ }^{-1}=-v/\eta (v)=-v/(1+\rho v)\text{ for }x\in \mathbb{G}_{\rho }%
\text{ (as }v\neq \rho ^{\ast }).
\]%
Isomorphic maps of $\mathbb{G}$ are provided for $\rho =0$ by $\iota
:x\mapsto x$ onto $(\mathbb{R},+),$ and for $\rho >0$ by $\eta :x\rightarrow
1+\rho x$ onto $(\mathbb{R}_{+},\cdot ),$ since
\[
\eta (u)\eta (v)=(1+\rho u)(1+\rho v)=1+\rho \lbrack u+v(1+\rho u)]=\eta
(u\circ _{\eta }v).
\]%
The rest follows since $\rho >0$ and $x>-1/\rho $ imply $\eta (x)>0.$ Also,
as $\rho \rho ^{\ast }=-1,$%
\begin{eqnarray*}
(2\rho ^{\ast }-s)+(2\rho ^{\ast }-t)(1+\rho (2\rho ^{\ast }-s)) &=&4\rho
^{\ast }-s-t+(2\rho ^{\ast }-t)(-2-\rho s)) \\
&=&s+t(1+\rho s)=s\circ _{\rho }t \\
&=&\pi ^{2}(s\circ _{\rho }t),
\end{eqnarray*}%
(as $\pi ^{2}=\iota $) and $|\eta _{\ast }^{\rho }(\pi (t))|=|1+\rho \lbrack
-t-2/\rho ]|=|-1-\rho t|=\eta _{\ast }^{\rho }(t),$ for $t\in (\rho ^{\ast
},\infty ).$ $\square $

\bigskip

\noindent \textbf{Remarks}. 1. For $\rho \neq 0$, $\mathbb{G}^{\rho }$ is
typified (rescaling its domain$)$ by the case $\rho =1,$ where%
\[
a\circ _{1}b=(1+a)(1+b)-1:\qquad (\mathbb{G}^{1},\circ _{1})=(\mathbb{R}%
^{\ast },\cdot )-1,
\]%
and the isomorphism is a shift (cf. [Pop, \S 3]), i.e. the groups are
conjugate. This is the classical \textit{circle group }above.

\noindent 2. For $\rho >0,$ note that $u\in \mathbb{G}_{+}^{\rho }\cap
(0,\infty )$ has $u_{\circ }^{-1}\in \mathbb{G}_{+}^{\rho }\cap (-1/\rho
,0). $

\noindent 3. Since $\eta (t_{\circ }^{-1})=1/\eta (t)$, $t_{\circ
}^{-1}\circ v=(v-t)/\eta (t),$ and so the convolution $t\ast v:=v\circ
t_{\circ }^{-1}$ is the asymptotic form of the Beurling convolution $%
(v-t)/\varphi (t)$ occurring in the Beurling Tauberian Theorem (\S 4) for $%
\varphi \in SN$.

\noindent 4. For $\rho >0,$ the inverse $\eta ^{-1}(y)=(y-1)/\rho $ maps $%
(0,\infty )$ onto $\mathbb{G};$ moreover, $\eta ^{-1}$ is \textit{%
super-additive} on $(1,\infty )$, i.e. for $x,y\geq 1=1_{\mathbb{R}^{\ast
}}, $
\[
\eta ^{-1}(x)+\eta ^{-1}(y)\leq \eta ^{-1}(xy),
\]%
as
\[
0\leq \rho ^{2}\eta ^{-1}(x)\eta ^{-1}(y)=(xy-1)-(x-1)-(y-1)=\rho \eta
^{-1}(xy)-\rho \eta ^{-1}(x)-\rho \eta ^{-1}(y);
\]%
it is also super-additive on $(0,1).$

\bigskip

Below we list further useful arithmetic facts including the iterates $%
a_{\varphi x}^{n+1}=a_{\varphi x}^{n}\circ _{\varphi x}a$ with $a_{\varphi
x}^{1}=a$. To avoid excessive bracketing, the usual arithmetic operations
below bind more strongly than Popa operations.

\bigskip

\noindent \textbf{Proposition 2 }(Arithmetic of Popa operations).%
\[
\begin{array}{cc}
\text{i)} & a_{\varphi x}^{0}=1_{\varphi x}=0;\qquad a\circ _{\varphi
x}a_{\varphi x}^{-1}=0\qquad \text{for }a_{\varphi x}^{-1}:=(-a)/\eta
_{x}^{\varphi }(a); \\
\text{ii)} & x\circ _{\varphi }(b\circ _{\varphi x}a)=y\circ _{\varphi
}a,\qquad \text{for }y:=x\circ _{\varphi }b; \\
\text{iii)} & x\circ _{\varphi }(b\circ _{\eta }a)=y\circ _{\varphi }a\eta
(b)/\eta _{x}^{\varphi }(b)\qquad \text{for }y:=x\circ _{\varphi }b; \\
\text{iv)} & x=y\circ _{\varphi }b_{\varphi x}^{-1}\qquad \text{for }%
y:=x\circ _{\varphi }b; \\
\text{v)} & \eta _{x}^{\varphi }(a_{\varphi
x}^{m})=\prod\limits_{k=0}^{m-1}\eta _{y_{k}}^{\varphi }(a),\qquad \text{for
the iterates }a_{\varphi x}^{n}\text{ and }y_{k}=x\circ _{\varphi
}a_{\varphi x}^{k},\text{ }(k=0,...,m-1).%
\end{array}%
\]%
\noindent \textit{Proof.} (i) Here $1_{\varphi x}$ denotes the neutral
element of the operation $\circ _{\varphi x},$ which is $0,$ since $\eta
_{x}^{\varphi }(0)=1$ (so that $0\circ _{\varphi x}t=t,$ while $s\circ
_{\varphi x}0=s$). So $a_{\varphi x}^{1}=a_{\varphi x}^{0}\circ _{\varphi
x}a.$%
\[
a\circ _{\varphi x}a_{\varphi x}^{-1}=a+a_{\varphi x}^{-1}\eta _{x}^{\varphi
}(a)=0.
\]

\noindent (ii) For $y=x\circ _{\varphi }b,$
\[
x\circ _{\varphi }(b\circ _{\varphi x}a)=x\circ _{\varphi }(b+a\eta
_{x}^{\varphi }(b))=x+b\varphi (x)+a\varphi (x+b\varphi (x))=y\circ
_{\varphi }a.
\]%
\noindent (iii) As in the preceding step for (ii),
\[
x+(b+a\eta (b))[\varphi (x)/\varphi (x\circ _{\varphi }b)]\varphi (x\circ
_{\varphi }b)=y\circ _{\varphi }a\eta (b)/\eta _{x}^{\varphi }(b).
\]

\noindent (iv) For $y=x\circ _{\varphi }b,$ using $b_{\varphi
x}^{-1}=-b/\eta _{x}^{\varphi }(b)$ from (i),%
\[
x=y-b\varphi (x)=y-[b\varphi (x)/\varphi (y)]\varphi (y)=y\circ _{\varphi
}b_{\varphi x}^{-1}.
\]

\noindent (v) For $m=1$ both sides agree since by (i) $y_{0}=x.$ Proceed by
induction, using (ii):%
\begin{eqnarray*}
\eta _{x}^{\varphi }(a_{\varphi x}^{m+1}) &=&\varphi (x\circ _{\varphi
}(a_{\varphi x}^{m}\circ _{\varphi x}a))/\varphi (x)=\varphi (y_{m}\circ
_{\varphi }a))/\varphi (x) \\
&=&[\varphi (y_{m}\circ _{\varphi }a)/\varphi (y_{m})]\varphi (x\circ
_{\varphi }a_{\varphi x}^{m})/\varphi (x)=\eta _{y_{m}}^{\varphi }(a)\eta
_{x}^{\varphi }(a_{\varphi x}^{m}).\text{ }\square
\end{eqnarray*}

\section{Extension to Beurling's Tauberian Theorem}

Theorem 2 below extends one proved by Beurling in lectures in 1957; see e.g.
[Kor, IV.11] for references. Bingham and Goldie [BinG2] extended Beurling's
result by replacing the Lebesgue integrator $H(y)dy$ below by a suitable
Lebesgue-Stieltjes integrator $dU(y),$ and demanding more of the Wiener
kernel (than just non-vanishing of its Fourier transform), and gave a
corollary for Beurling moving averages.

Here we extend the class of Beurling convolutions applied in the other term
of the integrand, replacing $\varphi \in BSV$ by $\varphi \in SE,$ so
widening the application to moving averages, as we note below. With the
following `Beurling notation' for Lebesgue and Stieltjes integrators,%
\begin{eqnarray*}
F\ast _{\varphi }H(x):= &&\int F\Bigl(\frac{x-u}{\varphi (x)}\Bigr)H(u)\frac{%
\mathrm{d}u}{\varphi (x)}=\int F(-t)H(x+t\varphi (x))\mathrm{d}t, \\
F\ast _{\varphi }dU(x):= &&\int F\Bigl(\frac{x-u}{\varphi (x)}\Bigr)\frac{%
\mathrm{d}U(u)}{\varphi (x)}=\int F(-t)dU(x+t\varphi (x))\mathrm{d}t,
\end{eqnarray*}%
reducing for $\varphi \equiv 1$ to their classical counterparts
\[
F\ast H(x)=\int F(x-t)H(t)\mathrm{d}t,\qquad \text{ }F\ast dU(x)=\int F(x-t)%
\mathrm{d}U(t),
\]%
we recall Wiener's theorem for the Lebesgue and the Lebesgue-Stieltjes
integrals. The latter uses the class $\mathcal{M}$ of continuous functions
(see Widder [Wid, V.12]; cf. [Wie, II.10]) with norm:%
\[
||f||:=\sup_{y\in \mathbb{R}}\sum_{n\in \mathbb{Z}}\sup_{x\in \lbrack
0,1]}|f(x+y+n)|<\infty ,
\]%
and places a uniform bounded-variation restriction on the integrator $U$ as
follows. Denote by $|\mu _{x}|$ the usual norm of the charge (signed
measure) generated from the function $y\mapsto U_{x}(x\circ _{\varphi
}y)/\varphi (x);$ then there should exist $\delta >0$ and $M<\infty $ with%
\begin{equation}
\sup_{x,y\in \mathbb{R}}|\mu _{x}|(I_{\delta }^{+}(y))\leq M,  \tag{$BV$}
\end{equation}%
where $I_{\delta }^{+}(y):=[y,y+\delta ).$ It will be convenient to refer to
the following conditions as $x\rightarrow \infty ,$ with or without the
subscript $\varphi $ (the latter when $\varphi \equiv 1):$%
\begin{equation}
K\ast _{\varphi }H(x)\rightarrow c\int K(y)\mathrm{d}y,\quad K\ast _{\varphi
}dU(x)\rightarrow c\int K(y)\mathrm{d}y.  \tag{$K\ast _{\varphi }H$/$U$}
\end{equation}

\noindent \textbf{Theorem W (Wiener's Tauberian Theorem)}.\textit{\ For }$%
K\in L_{1}(\mathbb{R})$\textit{\ }(\textit{resp. }$K\in \mathcal{M}$)\textit{%
\ with }$\hat{K}$\textit{\ non-zero on }$\mathbb{R}$\textit{:}

\textit{if }$H$\textit{\ is bounded }(\textit{resp. }$H\in \mathcal{M}$)%
\textit{, and }($K\ast H$)\textit{, resp. }($K\ast U$), \textit{holds, then
for all }$F\in L_{1}(\mathbb{R})$ (\textit{resp. }$F\in \mathcal{M}$)\textit{%
,}
\[
F\ast H(x),\text{ resp. }F\ast \mathrm{d}U(x)\rightarrow c\int F(t)\mathrm{d}%
t\qquad (x\rightarrow \infty ).
\]

\bigskip

\noindent \textbf{Theorem B (Beurling's Tauberian theorem)}.\textit{\ For }$%
K\in L_{1}(\mathbb{R})$\textit{\ with }$\hat{K}$\textit{\ non-zero on }$%
\mathbb{R}$\textit{, and }$\varphi $ \textit{Beurling slowly varying,}%
\begin{equation}
\varphi (x+t\varphi (x))/\varphi (x)\rightarrow 1,\qquad (x\rightarrow
\infty )\qquad (t\geq 0):  \tag{$BSV$}
\end{equation}

\textit{if }$H$\textit{\ is bounded, and }($K\ast _{\varphi }H$) \textit{%
holds, then for all }$F\in L_{1}(\mathbb{R})$
\[
F\ast _{\varphi }H(x)\rightarrow c\int F(y)\mathrm{d}y\qquad (x\rightarrow
\infty ).
\]%
We recommend the much later, slick, and elegant proof in [Kor, IV.11].

\bigskip

\noindent \textbf{Theorem BG1 (}$\mathbf{LS}$\textbf{-Extension to
Beurling's Tauberian theorem}, [BinG2, Th. 8]).\textit{\ If }$\varphi \in
BSV $\textit{, }$K\in \mathcal{M}$\textit{\ with }$\hat{K}$\textit{\
non-zero on }$\mathbb{R}$\textit{, }$U$\textit{\ satisfies }($BV$) \textit{%
and }($K\ast _{\varphi }U$) \textit{holds}

\textit{--- then for all }$G\in \mathcal{M}$\textit{,}
\[
G\ast _{\varphi }\mathrm{d}U(x)\rightarrow c\int G(y)\mathrm{d}y\qquad
(x\rightarrow \infty ).
\]%
We show how to amend the proof of Th. BG1\ in [BinG2] (similar in essence to
that cited above in [Kor, IV.11]) to obtain the following.

\bigskip

\noindent \textbf{Theorem 2 (Extension to Beurling's Tauberian theorem)}.%
\textit{\ If }$\varphi \in SE$, \textit{i.e. locally uniformly in }$t$%
\begin{equation}
\varphi (x+t\varphi (x))/\varphi (x)\rightarrow \eta (t)\in GS,\qquad
(x\rightarrow \infty )\qquad (t\geq 0),  \tag{$SE$}
\end{equation}%
\textit{\ }$K\in L_{1}(\mathbb{R})$ (\textit{resp. }$K\in \mathcal{M}$%
\textit{\ ) with }$\hat{K}$\textit{\ non-zero on }$\mathbb{R}$\textit{, }$H$%
\textit{\ is bounded }(\textit{resp.}$U$\textit{\ satisfies }($BV$)) \textit{%
and }($K\ast _{\varphi }H$), \textit{resp.} ($K\ast _{\varphi }U$),\textit{\
holds}

\textit{\ --- then for all }$G\in L_{1}(\mathbb{R})$ (\textit{resp. }$G\in
\mathcal{M}$)
\[
G\ast _{\varphi }H(x)\rightarrow c\int G(y)\mathrm{d}y,\qquad \text{resp. }%
G\ast _{\varphi }\mathrm{d}U(x)\rightarrow c\int G(y)\mathrm{d}y\qquad
(x\rightarrow \infty ).
\]

\bigskip

\noindent \textit{Proof.} We consider the Lebesgue-Stieltjes case (the
Lebesgue case is similar, but simpler). For fixed $a$ and with $K$ as in the
Theorem, set $K_{a}(s):=K(s-a),$ and take%
\[
t:=(s-a)/\eta _{x}(a),\text{ }dt=ds/\eta _{x}(a)\text{ and }s=a+t\eta
_{x}(a)=a\circ _{\varphi x}t.
\]%
Then for $y=x+a\varphi (x),$ by Prop. 2(ii), $x\circ _{\varphi }(a\circ
_{\varphi x}t)=y\circ _{\varphi }a$ and so%
\begin{eqnarray*}
K_{a}(s)U(x\circ _{\varphi }s) &=&K(t\eta _{x}(a))U(x\circ _{\varphi
}(a\circ _{\varphi x}t)) \\
&=&K(t\eta _{x}(a))U(y\circ _{\varphi }t).
\end{eqnarray*}%
So, as in [BinG2], for $K$ continuous ($K\in \mathcal{M}$), with $A$ as in
[BinG2] for $c$ above%
\begin{eqnarray*}
\int K_{a}(s)\mathrm{d}U(x\circ _{\varphi }s) &=&\eta _{x}(a)\int K(t\eta
_{x}(a))\mathrm{d}U(y\circ _{\varphi }t)\rightarrow A\int K(t\eta (a))\eta
(a)\mathrm{d}t \\
&=&A\int K(u)\mathrm{d}u,\text{ for }u:=t\eta (a).
\end{eqnarray*}%
Now continue with the proof verbatim as in [BinG2]. $\square $

\bigskip

\noindent \textbf{Corollary 3 }([BinG2, \S 5 Cor. 2] for $\varphi \in SN$).
\textit{For }$\varphi \in SE,$ \textit{if }$U$\textit{\ is non-decreasing
and for some }$\delta >0$
\[
\sup_{x,y\in \mathbb{R}}[U_{x}(x\circ _{\varphi }(y+\delta ))-U_{x}(x\circ
_{\varphi }y]/\varphi (x)<\infty
\]%
\textit{--- then }$(K\ast _{\varphi }U)$\textit{\ holds for some }$c$\textit{%
\ and Wiener kernel }$K\in \mathcal{M}$ \textit{iff for some }$c_{U}$\textit{%
\ either of the following holds:}%
\[
(\Delta _{t}^{\varphi }U/\varphi )(x)\equiv \lbrack U(x\circ _{\varphi
}t)-U(x)]/\varphi (x)\rightarrow c_{U}t\qquad (x\rightarrow \infty )\qquad
(t>0),
\]%
\[
(\Delta _{t}^{\varphi }U/\varphi )(x)\rightarrow c_{U}t\qquad (x\rightarrow
\infty )\qquad \text{for two incommensurable }t.
\]

\noindent \textit{Proof.}\textbf{\ }Repeat verbatim the proof in [BinG2, \S %
5 Cor. 2], using $H(x)=t^{-1}\mathbf{1}_{[0,t]}(x),$ with $\mathbf{1}%
_{[0,t]} $ the indicator function of the interval $[0,t].$ $\square $

\section{Uniformity, semicontinuity}

To motivate our results below of limsup convergence type, we use the
following weak notion of uniformity: say that $f_{n}\rightarrow f$ \textit{%
uniformly near} $t$ if for every $\varepsilon >0$ there is $\delta >0$ and $%
m\in \mathbb{N}$ such that
\[
f(t)-\varepsilon <f_{n}(s)<f(t)+\varepsilon \text{ for }n>m\text{ and }s\in
I_{\delta }(t),
\]%
where $I_{\delta }(t):=(t-\delta ,t+\delta ).$ For instance, if $\varphi \in
SN$, $x_{n}$ divergent, $f(s)\equiv 1,$ and $f_{n}(s):=\varphi
(x_{n}+s\varphi (x_{n}))/\varphi (x_{n}),$ then `$f_{n}\rightarrow f$
uniformly near $t$ for all $t>0.$'

The notion above is easier to satisfy than Hobson's `uniform convergence
\textit{at} $t$' which replaces $f(t)$ above by $f(s)$ twice, [Hob, p.110];
suffice it to refer to $f_{n}\equiv 0,$ and $f$ with $f(0)=0$ and $f\equiv 1$
elsewhere. (See also Klippert and Williams [KliW], where though Hobson's
condition is satisfied at all points of a set, the choice of $\delta $
cannot itself be uniform in $t.)$

Our notion of uniformity may be equivalently stated in limsup language, as
follows, bringing to the fore the underlying \textit{uniform upper and lower
semicontinuity}. The proof of the next result is routine, so we omit it here%
\footnote{%
See the Appendix for details.}; but its statement will be useful in the
development below.

\bigskip

\noindent \textbf{Proposition 3 (Uniform semicontinuity).} \textit{If }$%
f_{n}\rightarrow f$\textit{\ pointwise, then }$f_{n}\rightarrow f$\textit{\
converges uniformly near }$t$\textit{\ iff}
\begin{eqnarray*}
f(t) &=&\lim_{\delta \downarrow 0}\limsup_{n}\sup \{f_{n}(s):s\in I_{\delta
}(t)\} \\
&=&\lim_{\delta \downarrow 0}\liminf_{n}\inf \{f_{n}(s):s\in I_{\delta
}(t)\}.
\end{eqnarray*}%
Again putting $I_{\delta }^{+}(t):=[t,t+\delta )$, we may now consider the
\textit{one-sided limsup-sup condition} at $t:$%
\begin{equation}
f(t)=f_{+}(t)\text{ with }f_{+}(t):=\lim_{\delta \downarrow
0}\limsup_{n}\sup \{f_{n}(s):s\in I_{\delta }^{+}(t)\}.  \label{unif+}
\end{equation}

The next result is akin to the Dini/P\'{o}lya-Szeg\H{o} monotone convergence
theorems (respectively [Rud,7.13], for monotone convergence of continuous
functions to a continuous pointwise limit, and [PolS, Vol. 1 p.63, 225,
Problems II 126, 127], or Boas [Boa, \S 17, p. 104-5], when the functions
are monotone); here we start with one-sided assumptions on the domain and
range, and conclude via a category argument by improving to a two-sided
condition.

\bigskip

\noindent \textbf{Proposition 4 (Uniform Upper semicontinuity).} \textit{If
quasi everywhere }$f_{n}$ \textit{converges pointwise to an upper
semicontinuous limit }$f$\textit{\ satisfying the one-sided condition }(\ref%
{unif+})\textit{\ quasi everywhere in its domain, then quasi everywhere }$f$%
\textit{\ is uniformly upper semicontinuous:}
\[
f(t)=\lim_{\delta \downarrow 0}\limsup_{n}\sup \{f_{n}(s):s\in I_{\delta
}(t)\}.
\]

\noindent \textit{Proof.} Take $\{I_{n}\}_{n\in \mathbb{N}}$ with $I_{0}=%
\mathbb{R}$ to be a sequence of open intervals that form a base for the
usual open sets of $\mathbb{R}.$ Let $\mathbb{D}$ be a countable dense
subset of the (co-meagre) intersection $S$ of the set on which $f_{n}$
converges pointwise and the set on which (\ref{unif+}) holds. Put%
\[
G^{k}(\varepsilon ):=\bigcup\nolimits_{m\in \mathbb{N}}\{(q,q+\delta ):q\in
\mathbb{D}\text{, }(q,q+\delta )\subseteq I_{k},\text{ }(\forall
n>m)(\forall s\in I_{\delta }^{+}(q))
\]%
\[
\lbrack f_{n}(s)<f(q)+\varepsilon ]\},
\]%
which is open. It is also dense in $I_{k}$: from any open interval $%
I\subseteq I_{k}$ choose $q\in \mathbb{D\cap }I$; as $q\in S\cap I,$ there
exist $N_{k}\in \mathbb{N}$ and $\delta >0$ such that $I_{\delta
}^{+}(q)\subseteq I$ and
\[
f_{n}(s)<f(q)+\varepsilon \text{\qquad (}n>N_{q},\text{ }s\in I_{\delta
}^{+}(q)\text{)};
\]%
so $(q,q+\delta )\subseteq I\cap G^{k}(\varepsilon ),$ i.e. $%
G^{k}(\varepsilon )$ meets $I.$ Consider $T_{k}:=\bigcap\nolimits_{%
\varepsilon \in \mathbb{Q}_{+}}G^{k}(\varepsilon )\subseteq I_{k}$; then, by
Baire's Theorem, $I_{k}\backslash T_{k}$ is meagre. Put%
\[
T:=T_{0}\backslash \bigcup\nolimits_{k\in \mathbb{N}}(I_{k}\backslash
T_{k}):\qquad T\cap I_{k}\subseteq T_{k}\qquad (k\in \mathbb{N}).
\]%
As $T$ is co-meagre, we may assume w.l.o.g. that the one-sided uniformity
condition (\ref{unif+}) holds on $T$.

Given $\varepsilon >0$ and $t\in T,$ by upper-semicontinuity of $f$ at $t,$
pick $r\in \mathbb{N}$ such that $t\in I_{r}$ and $f(u)<f(t)+\varepsilon $
for all $u\in I_{r}.$ Now, as $t\in T_{r}$, $t\in G^{r}(\varepsilon ),$ so
we may pick $q\in \mathbb{D}\cap I_{r}$ and $\delta >0$ with $t\in
(q,q+\delta )\subseteq I_{r}$ and $m\in \mathbb{N}$ such that%
\[
f_{n}(s)<f(q)+\varepsilon \text{\qquad }(n>m,\text{ }s\in I_{\delta
}^{+}(q)),
\]%
again as $q\in S.$ Now choose $d>0$ such that $I_{d}(t)\subseteq (q,q+\delta
).$ Then for $n>m$ and $s\in I_{d}(t)$
\[
f_{n}(s)<f(q)+\varepsilon <f(t)+2\varepsilon ,
\]%
since $q\in I_{r}.$ As $\varepsilon >0$ was arbitrary,%
\[
f(t)=\lim_{\delta \downarrow 0}\limsup_{n}\sup \{f_{n}(s):s\in I_{\delta
}(t)\}\text{ for }t\in T.\text{\qquad }\square
\]

\bigskip

Before proceeding further we need to extend the Beurling function $\eta
^{\varphi }$ some way to the left of the (natural) origin as follows
(recalling from \S 1 the condition $(SE_{\mathbb{A}})$); cf. BGT (2.11.2).
Here we see the critical role of the \textit{Popa} origin $\rho ^{\ast
}=-\rho ^{-1}$ of \S \S 1,3: the domain of the limit operation $%
\lim_{x\rightarrow \infty }\eta _{x}^{\varphi }(s),$ used to extend $\eta
^{\varphi },$ is $\mathbb{G}_{+}^{\rho }$, i.e. $s$ has to be to the right
of the Popa origin.

\bigskip

\noindent \textbf{Lemma 1 (Uniform Involutive Extension).} \textit{For }$%
\varphi \in SE$,\textit{\ }$\circ =\circ _{\rho }$ \textit{with }$\rho =\rho
_{\varphi }>0,$ \textit{put}%
\[
\eta ^{\varphi }(t_{\circ }^{-1})=\eta ^{\varphi }(-t/\eta ^{\varphi
}(t)):=1/\eta ^{\varphi }(t),\qquad (t>0);
\]%
\textit{then }$(SE_{\mathbb{A}})$ \textit{holds for} $\mathbb{A=G}_{+}^{\rho
}\mathbb{=}(\rho ^{\ast },\infty ).$ \textit{Moreover, this is a maximal
positive extension: for each} $s<\rho ^{\ast },$ \textit{assuming }$\varphi
(x+s\varphi (x))>0$\textit{\ is defined for all large }$x$\textit{,} \textit{%
\ }%
\[
\lim_{x\rightarrow \infty }\eta _{x}^{\varphi }(s)=\lim_{x\rightarrow \infty
}\varphi (x+s\varphi (x))/\varphi (x)=0=\eta (\rho ^{\ast }).
\]

\bigskip

\noindent \textit{Proof.} Fixing $t>0$ and taking $y:=x+t\varphi (x)$ and $%
s_{x}=t/\eta _{x}(y)$ gives $x=y-s_{x}\varphi (y).$ Now%
\[
1/\eta _{x}(y)=\frac{\varphi (x)}{\varphi (x+t\varphi (x))}\rightarrow
1/\eta ^{\varphi }(t)=\eta ^{\varphi }(-t/\eta (t)).
\]%
So $s_{x}\rightarrow s=t/\eta (t)$ and%
\[
\frac{\varphi (y-s_{x}\varphi (y))}{\varphi (y)}=\frac{\varphi (x)}{\varphi
(x+t\varphi (x))}\rightarrow 1/\eta ^{\varphi }(t)=\eta ^{\varphi }(-t/\eta
^{\varphi }(t))=\eta ^{\varphi }(-s).
\]%
So for $s>0$ with $\eta ^{\varphi }(-s)>0$ and $y$ so large that $%
y(1-s\varphi (y)/y)>0,$%
\[
\frac{\varphi (y-s\varphi (y))}{\varphi (y)}\rightarrow \eta (-s)\text{
locally uniformly in }s\text{ for }\eta (-s)>0.\text{ }
\]

As for the maximality assertion (even allowing $\mathbb{R}$ to be the domain
of $\varphi $), since%
\[
0\leq \liminf_{x\rightarrow \infty }\varphi (x+s\varphi (x))/\varphi (x)\leq
\limsup_{x\rightarrow \infty }\varphi (x+s\varphi (x))/\varphi (x),
\]%
suppose there are $s<\rho ^{\ast }$ and a divergent sequence $x_{n}$ with%
\[
\tilde{\eta}(s):=\lim_{n\rightarrow \infty }\varphi (x_{n}+s\varphi
(x_{n}))/\varphi (x_{n})>0,
\]%
including here the case $\tilde{\eta}(s)=+\infty .$ As above, take $%
y_{n}=x_{n}+s\varphi (x_{n})$ and $s_{n}=-s\varphi (x_{n})/\varphi
(x_{n}+s\varphi (x_{n}))>0;$ then $x_{n}=y_{n}-s\varphi
(x_{n})=y_{n}+s_{n}\varphi (y_{n})$ and $s_{n}\rightarrow -s\tilde{\eta}%
(s)^{-1}\geq 0,$ strictly so unless $\tilde{\eta}(s)=+\infty .$ So
\[
\tilde{\eta}(s)^{-1}=\lim \varphi (y_{n}+s_{n}\varphi (y_{n}))/\varphi
(y_{n})=\eta (-s\tilde{\eta}(s)^{-1})=1-\rho s\tilde{\eta}(s)^{-1}.
\]%
If $\tilde{\eta}(s)=+\infty ,$ this is already a contradiction. If $0<\tilde{%
\eta}(s)<\infty $ cross-multiplying by $\tilde{\eta}(s),$ yields $\tilde{\eta%
}(s)=1+\rho s<1+\rho \rho ^{\ast }=0,$ again a contradiction. $\square $

\bigskip

\noindent \textbf{Remark. }For $s>0$ and large enough $y$ the expression $%
y-s\varphi (y)$ is positive provided $s<\lim \inf x/\varphi (x),$ that is
for $s>\rho ^{\ast }.$ This corresponds to $\varphi (x)=O(x);$ if, however,
as in BGT \S 2.11, $\varphi (x)=o(x),$ then $\rho _{\varphi }=0,$ so that $%
\rho ^{\ast }=-\infty ,$ and so $s$ may be arbitrary.

\bigskip

\noindent \textbf{Definitions. }Recalling (\S 1) that $\Delta _{t}^{\varphi
}h(x):=h(x+t\varphi (x))-h(x),$ and, taking limits here and below as $%
x\rightarrow \infty $ (rather than sequentially as $n\rightarrow \infty ),$
put for $\varphi \in SE$ and $\rho =\rho _{\varphi }$%
\begin{eqnarray*}
\mathbb{A}^{\varphi }:= &\{t>&\rho ^{\ast }:\Delta _{t}^{\varphi }h\text{
converges to a finite limit}\}, \\
\mathbb{A}_{\text{u}}:= &\{t>&\rho ^{\ast }:\Delta _{t}^{\varphi }h\text{
converges to a finite limit locally uniformly near }t\}.
\end{eqnarray*}%
(For $\mathbb{A}^{\varphi }\subseteq \mathbb{G}_{+}^{\rho }$, see Lemma 1
above and Prop. 6 below.) So $0\in \mathbb{A}^{\varphi },$ but we cannot yet
assume either that $\mathbb{A}^{\varphi }$ is a subgroup, or that $0\in
\mathbb{A}_{\text{u}}$, a critical point in Proposition 6 below. In the
Karamata case $\varphi \equiv 1,$ $\mathbb{A}^{\varphi }=\mathbb{A}^{1}$ is
indeed a subgroup (see [BinO12, Prop. 1]).

For $t\in \mathbb{A}^{\varphi }$ put%
\begin{equation}
K(t):=\lim_{x\rightarrow \infty }\Delta _{t}^{\varphi }h.  \tag{$K$}
\end{equation}%
So $K(0)=0.$

Proposition 5 below is included to help in reading the subsequent
Proposition 6 -- dedicated to checking when $\mathbb{A\subseteq G}$ is a
subgroup of a Popa group -- which needs a sequential characterization of
uniform convergence near a non-zero $t$ (as $t_{n}\rightarrow t$ iff $%
c_{n}=t_{n}/t\rightarrow 1$); the proof is routine, so omitted.

\bigskip

\noindent \textbf{Proposition 5.} $h(x+t\varphi (x))-h(x)$ \textit{converges
locally (right-sidedly) uniformly} \textit{to }$K(t)$\textit{\ near }$t\neq
0 $\textit{, iff for each divergent }$x_{n}$\textit{\ and any }$%
c_{n}\rightarrow 1$\textit{\ }($c_{n}\downarrow 1)$\textit{\ }%
\[
h(x_{n}\circ _{\varphi }c_{n}t)-h(x_{n})\rightarrow K(t);
\]%
\textit{then, taking suprema over sequences }$c=\{c_{n}\}\downarrow 1$%
\textit{\ and }$x=\{x_{n}\}\rightarrow \infty $,\textit{\ }%
\[
K(t):=\sup_{c,x}\{\limsup_{n\rightarrow \infty }h(x_{n}\circ _{\varphi
}c_{n}t)-h(x_{n})\}.
\]

\bigskip

\noindent \textbf{Proposition 6.} \textit{For }$\varphi \in SE,$ $\mathbb{A}%
_{\text{u}}$ \textit{is a subgroup of }$\mathbb{G}_{+}^{\rho }$ \textit{for }%
$\rho =\rho _{\varphi }$ \textit{iff} $0\in \mathbb{A}_{\text{u}}$; \textit{%
then }$K:(\mathbb{A}_{\text{u}},\circ )\rightarrow (\mathbb{R},+),$\textit{\
defined by }$(K)$\textit{\ above, is a homomorphism}.

\bigskip

\noindent \textit{Proof.} We show that $v\circ _{\eta }u\in \mathbb{A}_{%
\text{u}}$ for $u,v\in \mathbb{A}_{\text{u}}$ with $v\circ _{\eta }u\neq 0,$
and that $\mathbb{A}_{\text{u}}$ is closed under inverses $u_{\circ }^{-1}$
for non-zero $u$, so it is a subgroup of $\mathbb{G}$ iff $1_{\mathbb{G}%
}=0\in \mathbb{A}_{\text{u}}$. For $u,v$ $\in \mathbb{A}_{\text{u}},$ since $%
\eta _{x}(v)=\varphi (x+v\varphi (x))/\varphi (x)\rightarrow \eta (v),$%
\[
u_{v}:=u\eta (v)/\eta _{x}(v)\rightarrow u,
\]%
and so with $y=x\circ _{\varphi }v,$ since by Prop. 2(iii) $x\circ _{\varphi
}(v\circ _{\eta }u)=y\circ _{\varphi }u_{v},$%
\begin{eqnarray*}
h(x\circ _{\varphi }(v\circ _{\eta }u))-h(x) &=&[h(y\circ _{\varphi
}u_{v})-h(y)]+[h(x\circ _{\varphi }v)-h(x)] \\
&\rightarrow &K(u)+K(v),
\end{eqnarray*}%
i.e.%
\[
K(v\circ _{\eta }u)=\lim [h(x\circ _{\varphi }(v\circ _{\eta
}u))-h(x)]=K(u)+K(v).
\]%
As the convergence at $u,v$ on the right occurs uniformly near $u,v$
respectively, this is uniform near $v\circ u,$ using Prop. 5 provided $%
v\circ u\neq 0$.

For non-zero $t$ $\in \mathbb{A}_{\text{u}},$ this time put $y:=x\circ
_{\varphi }t;$ then, by Prop. 2(iv)$,$ $x=y\circ _{\varphi }t_{\varphi
x}^{-1},$ so%
\[
h(y\circ _{\varphi }t_{\varphi x}^{-1})-h(y)=[h(x)-h(y)]=-[h(x\circ
_{\varphi }t)-h(x)]\rightarrow -K(t).
\]%
So, since $t_{\varphi x}^{-1}=-t/\eta _{x}(t)\rightarrow -t/\eta (t),$%
\[
K(t_{\circ }^{-1})=K(-t/\eta (t))=\lim [h(y\circ _{\varphi }t_{\varphi
x}^{-1})-h(y)]=-K(t).
\]%
That is $t_{\circ }^{-1}\in \mathbb{A}_{\text{u}}$ (and $K(t_{\circ
}^{-1})=-K(t));$ again this is locally uniform at $t\neq 0,$ using Prop. 5. $%
\square $

\bigskip

Theorem 3 (a corollary of Proposition 6) and Theorem 4 below, together with
the results of \S 6 below, are of \textit{dichotomy }type.\textit{\ }The
theme is that uniformity holds nowhere or (under assumptions) everywhere.

\bigskip

\noindent \textbf{Theorem 3.} \textit{If} $\mathbb{A}_{\text{u}}$ \textit{is
non-empty, then} $0\in \mathbb{A}_{\text{u}}$ \textit{and so} $\mathbb{A}_{%
\text{u}}$ \textit{is a subgroup. }

\noindent \textit{In particular, for }$h(t)=\log \varphi (t)$\textit{, if }$%
\eta _{x}^{\varphi }(t)\rightarrow \eta (t)$ \textit{locally uniformly near }%
$t$\textit{\ for some }$t>0$\textit{, then this convergence is locally
uniform near }$t$\textit{\ for all }$t\geq 0.$

\bigskip

\noindent \textit{Proof.} Choose $s\in \mathbb{A}_{\text{u}}$, which without
loss of generality is non-zero (otherwise there is nothing to prove). So, as
above, $t:=-s/\eta (s)\in \mathbb{A}_{\text{u}}$. For arbitrary $%
z_{n}\rightarrow 0$ and $x_{n}$ divergent, take $s_{n}:=s+z_{n}$ and $%
t_{n}:=-s/\eta _{x(n)}(s_{n})\rightarrow t;$ then $y_{n}=x_{n}+(s+z_{n})%
\varphi (x_{n})$ is divergent. So (since $s\varphi (x_{n})=(s/\eta
_{x(n)}(s_{n}))\varphi (y_{n}))$%
\[
h(x_{n}+z_{n}\varphi (x_{n}))-h(x_{n})=h(x_{n}+(s+z_{n})\varphi
(x_{n})-(s/\eta _{x(n)}(s_{n}))\varphi (y_{n}))-h(x_{n}),
\]%
which (as $y_{n}=x_{n}+(s+z_{n})\varphi (x_{n}))$ is
\begin{eqnarray*}
&=&h(y_{n}-s/\eta _{x(n)}(s_{n})\varphi
(y_{n}))-h(y_{n})+h(x_{n}+s_{n}\varphi (x_{n}))-h(x_{n}) \\
&=&h(y_{n}\circ _{\varphi }t_{n})-h(y_{n})+h(x_{n}\circ _{\varphi }s_{n}) \\
&\rightarrow &h(t)+h(s)=h(s_{\varphi }^{-1})+h(s)=0.
\end{eqnarray*}%
So $\Delta _{t}^{\varphi }h$ converges locally near $0,$ i.e. $0\in \mathbb{A%
}_{\text{u}}$ -- a subgroup, by Prop. 6.

In particular, for $h=\log \varphi ,$%
\[
h(x_{n}+z_{n}\varphi (x_{n}))-h(x_{n})\rightarrow 0\text{ iff }\varphi
(x_{n}+z_{n}\varphi (x_{n}))/\varphi (x_{n})\rightarrow 1,
\]%
and $\mathbb{A}_{\text{u}}$ is non-empty as $\varphi \in SE.$ $\square $

\bigskip

The following result extends the Uniformity Lemma of [BinO10, Lemma 3].
Although the proof parallels the original, the current one-sided context
demands the closer scrutiny offered here. To describe more accurately the
convergence in $(K)$ above, we write%
\begin{eqnarray}
\Delta _{t}^{\varphi }h(x) &\rightarrow &K_{+}(t)\text{ if uniform near }t%
\text{ on the right,}  \TCItag{$K_{+}$} \\
\Delta _{t}^{\varphi }h(x) &\rightarrow &K_{-}(t)\text{ if uniform near }t%
\text{ on the left,}  \TCItag{$K_{-}$} \\
\Delta _{t}^{\varphi }h(x) &\rightarrow &K_{\pm }(t)\text{ if uniform near }t%
\text{.}  \TCItag{$K_{\pm }$}
\end{eqnarray}

\bigskip

\noindent \textbf{Lemma 2.}\textit{\ }(i) \textit{For }$\varphi \in SE:$

\noindent (a)\textit{\ if the convergence in }$(K)$ \textit{is uniform
(resp. right-sidedly uniform) near }$t=0$\textit{, then it is uniform (resp.
right-sidedly uniform) everywhere in }$\mathbb{A}^{\varphi }$ \textit{and
for }$u\in \mathbb{A}^{\varphi }\cap (0,\infty )$\textit{\ }%
\[
K_{+}(u)=K(u)+K_{+}(0);
\]

\noindent (b) \textit{if the convergence in }$(K)$\textit{\ is uniform near }%
$t=u\in \mathbb{A}^{\varphi }\cap (\mathbb{A}^{\varphi })_{\circ }^{-1}\cap
(0,\infty ),$\textit{\ then it is uniform near }$t=0:$%
\[
K_{\pm }(0)=K_{\pm }(u)+K(u_{\circ }^{-1});
\]

(ii) \textit{if }$\rho =0$ \textit{and }$\varphi \in SN$ \textit{is
monotonic increasing, and the convergence in }$(K)$\textit{\ is
right-sidedly uniform near }$t=u\in \mathbb{A}^{\varphi }\cap (0,\infty ),$%
\textit{\ then it is right-sidedly uniform near }$t=0:$%
\[
K_{+}(0)=K_{+}(u)+K(u_{\circ }^{-1}).
\]

\bigskip

\noindent \textit{Proof.} (i) (a) Suppose $(K)$ holds locally right-sidedly
uniformly (uniformly) near $t=0.$ Let $u\in \mathbb{A}^{\varphi }$ and $%
z_{n}\downarrow 0$ (resp. $z_{n}\rightarrow 0).$ For $x_{n}$ divergent ($%
x_{n}\rightarrow \infty $)$,$ $y_{n}:=x_{n}\circ _{\varphi
}u=x_{n}(1+u\varphi (x_{n})/x_{n})$ is divergent and
\begin{equation}
h(x_{n}\circ _{\varphi }(u+z_{n}))-h(x_{n})=h(x_{n}\circ _{\varphi
}u)-h(x_{n})+h(y_{n}\circ _{\varphi }z_{n}/\eta _{x(n)}(u))-h(y_{n}).
\tag{$\ast $}
\end{equation}%
Without loss of generality $\eta _{x(n)}(u)>0$ (all $n),$ since $u\in
\mathbb{G}_{+}^{\rho }$ and so%
\[
\eta _{x(n)}(u)\rightarrow \eta (u)>0;
\]%
then $z_{n}/\eta _{x(n)}(u)\downarrow 0$ (resp. $z_{n}/\eta
_{x(n)}(u)\rightarrow 0$). From $h(x_{n}\circ _{\varphi
}u)-h(x_{n})\rightarrow K(u),$ and the assumed uniform behaviour at the
origin, there is right-sidedly uniform (uniform) behaviour near $u.$ The
second statement follows on specializing to $u\in \mathbb{A}^{\varphi }\cap
(0,\infty )$ and taking limits in $(\ast )$.

(b) For the converse we argue as in Theorem 3. Suppose uniformity holds near
$u\in \mathbb{A}^{\varphi }\cap (\mathbb{A}^{\varphi })^{-1}\cap (0,\infty
); $ then $v:=u_{\circ }^{-1}=-u/\eta (u)\in \mathbb{A}^{\varphi }\cap (\rho
^{\ast },0).$ Let $z_{n}\rightarrow 0;$ then $z_{n}^{\prime }:=z_{n}/\eta
_{x(n)}(v)\rightarrow 0,$ as $\eta _{x(n)}(v)\rightarrow \eta (v).$ Also $%
(-v)/\eta _{x(n)}(v)\rightarrow (-v)/\eta (v)=v_{\circ }^{-1}=u,$ so%
\[
\lim (-v+z_{n})/\eta _{x(n)}(v)=u+0.
\]%
Taking $y_{n}:=x_{n}\circ _{\varphi }v$ ($<x_{n}$ for $v<0,$as here)%
\[
x_{n}\circ _{\varphi }z_{n}=(x_{n}\circ _{\varphi }v)\circ _{\varphi
}(-v+z_{n})/\eta _{x(n)}(v),
\]%
and%
\begin{eqnarray*}
h(x_{n}\circ _{\varphi }z_{n})-h(x_{n}) &=&h(y_{n}\circ _{\varphi
}(-v+z_{n})/\eta _{x(n)}(v))-h(y_{n})+h(x_{n}\circ _{\varphi }v)-h(x_{n}) \\
&\rightarrow &K(u)+K(v)=K(u)+K(u_{\circ }^{-1}),
\end{eqnarray*}%
where the convergence on the right is uniform in the first term and
pointwise in the second term.

(ii) When $\varphi \in SN$ is monotone, the argument in (b) above may be
amended to deal with right-sided convergence, as $1/\eta _{x(n)}(v)=\varphi
(x_{n})/\varphi (y_{n})\geq 1$ (for $v<0),$ and so $1/\eta _{x(n)}(v)$ tends
to $1$ from above, as $\rho =0.$ Also $z_{n}^{\prime }=z_{n},$ so if $%
z_{n}\downarrow 0,$ then $z_{n}\varphi (x_{n})/\varphi (y_{n})$ tends to $0$
from above, since $z_{n}\geq 0$ and
\[
(-v+z_{n})/\eta _{x(n)}(v)\text{ tends to }u\text{ from above,}
\]%
as $(-v)/\eta _{x(n)}(v)$ tends from above to $(-v)=u>0.$ From here the
argument is valid when `uniform' is replaced by `right-sidedly uniform'. $%
\square $

\bigskip

\noindent \textbf{Remark. }For $\varphi \in SE$ and $\eta =\eta ^{\varphi }$
write $\varphi \in SE^{+}/SE^{-}$ (for $u>0)$ respectively according as
\[
\varphi (x+u\varphi (x))/\varphi (x)\text{ tends to }\eta (u)\text{ from
below, or from above }
\]%
as $x\rightarrow \infty ,$ and likewise for $\varphi \in SN$ (with $\eta
^{\varphi }\equiv 1)$ and $SN^{-}.$ So if $\varphi \in SN$ and $\varphi $ is
increasing, then $\varphi \in SN^{-}$, since $\varphi (x+u\varphi
(x))>\varphi (x)$ for $u>0,$ so%
\[
\varphi (x+u\varphi (x))/\varphi (x)\text{ tends to }1\text{ from above.}
\]%
This was used in (ii) above, and extends to $SE.$ Of course $\eta \in
SE^{+}\cap SE^{-}.$

\bigskip

The next result leads from a one-sided condition to a two-sided conclusion.
This is the prototype of further such results, useful later.

\bigskip

\noindent \textbf{Theorem 4.} \textit{If the pointwise convergence }$(K)$%
\textit{\ holds on a co-meagre set in }$\mathbb{G}_{+}^{\rho }$ \textit{with
the limit function }$K$\textit{\ upper semicontinuous also on a co-meagre
set, and the one-sided condition}%
\begin{equation}
K(t)=\lim_{\delta \downarrow 0}\limsup_{x\rightarrow \infty }\sup
\{h(x+s\varphi (x))-h(x):s\in I_{\delta }^{+}(t)\}  \tag{$UNIF^{+}$}
\end{equation}%
\textit{\ holds at the origin -- then two-sided limsup convergence holds
everywhere: }%
\[
\mathbb{A}^{\varphi }=\mathbb{A}_{\text{u}}=\mathbb{G}_{+}^{\rho }.
\]

\noindent \textit{Proof.} The pointwise convergence assumption says $\mathbb{%
A}^{\varphi }$ is co-meagre (in \textit{\ }$\mathbb{G}_{+}^{\rho }$);
w.l.o.g. $\mathbb{A}^{\varphi }=(\mathbb{A}^{\varphi })_{\circ }^{-1},$
otherwise work below with the co-meagre set $\mathbb{A}^{\varphi }\cap (%
\mathbb{A}^{\varphi })_{\circ }^{-1}.$ Take $f(t):=K(t)$; then $%
f_{n}(t):=h(x_{n}\circ _{\varphi }t)-h(x_{n})\rightarrow f(t)$ holds
pointwise quasi everywhere on $\mathbb{A}^{\varphi }$. Since $(UNIF^{+})$
holds at $t=0,$ by Lemma 2(i)(a), it holds everywhere in $\mathbb{A}%
^{\varphi }$ and so quasi everywhere. By Proposition 4, its two-sided limsup
version holds quasi everywhere, and so at some point $u\in \mathbb{A}%
^{\varphi }\cap (\mathbb{A}^{\varphi })_{\circ }^{-1}\cap (0,\infty ).$ Then
by Lemma 2(i)(b) the two-sided limsup version holds at $0$, and so by Lemma
2(i)(a) it holds everywhere in $\mathbb{A}^{\varphi }.$ It now follows that $%
0\in \mathbb{A}^{\varphi }=\mathbb{A}_{\text{u}}$ and so $\mathbb{A}_{\text{u%
}}$ is a co-meagre subgroup of $\mathbb{G}_{+}^{\rho }$; so, by the
Steinhaus Subgroup Theorem (see [BinO9]), which applies here by Prop. 6, $%
\mathbb{A}_{\text{u}}=\mathbb{G}_{+}^{\rho }$. $\square $

\section{Dichotomy}

We continue with the setting of \S 5, but here we assume less about $\mathbb{%
A}^{\varphi }$ -- in place of being co-meagre we ask that it contains a
non-meagre Baire subset $S\subseteq $\textit{\ }$\mathbb{G}_{+}^{\rho }$.
This is a local version of the situation in \S 5 in that

\noindent (i) $S$ is locally co-meagre quasi everywhere, and

\noindent (ii) $\mathbb{A}^{\varphi }$ is non-meagre and contains a Baire
subset to witness this.

\noindent For general $h$ and $\varphi $ we cannot assume this happens.
However, under certain axioms of set-theory this will be guaranteed: see \S %
11. Now $\langle S\rangle ,$ the $\circ $-subgroup generated by $S,$ will of
course be $\mathbb{G}_{+}^{\rho }$, again by the Steinhaus Subgroup Theorem,
as in Theorem 4. So our aim here is to verify that $\mathbb{A}^{\varphi }$
is a subgroup by checking that $\mathbb{\mathbb{G}_{+}^{\rho }=}\langle
S\rangle \subseteq \mathbb{A}_{\text{u}}\subseteq \mathbb{A}^{\varphi }$.

\bigskip

\noindent \textbf{Theorem 5.} \textit{For }$\varphi \in SE$\textit{\ and }$h$%
\textit{\ Baire, if} $\mathbb{A}^{\varphi }$ \textit{contains a non-meagre
Baire subset, then} $\mathbb{A}^{\varphi }=\mathbb{G}_{+}$ \textit{and }$K$
\textit{is a homomorphism: }$K(u)=c\log (1+\rho t),$ \textit{for some }$c\in
\mathbb{R}$\textit{, }$(u\in \mathbb{G}_{+}),$ \textit{if }$\rho =\rho
_{\varphi }>0.$

\bigskip

Given our opening remarks, this reads as an extension of the Fr\'{e}%
chet-Banach Theorem on the continuity of Baire/measurable additive functions
-- for background see [BinO9]. The proof (see below) parallels Prop. 1 of
[BinO12], extending the cited result from the Karamata to the Beurling
setting, but now we need the Baire property to employ uniformity arguments
here.

Proposition 7 extends Theorem 7 (UCT) of [BinO10] and is crucial here.

\bigskip

\noindent \textbf{Proposition 7 (Uniformity).}\textit{\ Suppose }$S\subseteq
\mathbb{A}^{\varphi }$\textit{\ for some Baire non-meagre} $S$\textit{. Then
for Baire }$h$\textit{\ the convergence in}$\mathit{\ }$($K$) \textit{of \S %
5 is uniform near }$u=0$\textit{\ and so also near }$u=t$\textit{\ for }$%
t\in S,$ \textit{i.e. }$S\subseteq \mathbb{A}_{\text{u}}\subseteq \mathbb{A}%
^{\varphi }.$

\bigskip

\noindent \textit{Proof .} For each $n,$ define for $t>\rho ^{\ast }$ the
function $k_{n}(t):=h(n\circ _{\varphi }t)-h(n),$ which is Baire; then for $%
t\in \mathbb{A}^{\varphi }$
\[
K(t)=\lim\nolimits_{n}k_{n}(t),
\]%
and so $k=K|S$ is a Baire function with non-meagre domain. Now apply the
argument of Theorem 7 of [BinO10] to $S$ and $k$ as defined here (so that
Baire's Continuity Theorem [Oxt, Th. 8.1] applies to the Baire function $k$%
), giving uniform convergence near $u=0,$ so uniform convergence near any $%
u\in S,$ by Lemma 2(i)(a)$.$ $\square $

\bigskip

\noindent \textbf{Corollary 4.} \textit{If }$S\subseteq \mathbb{A}^{\varphi
} $ \textit{with }$S\subseteq \mathbb{G}_{+}^{\rho }$\textit{\ Baire and
non-meagre, and }$\rho \geq 0$\textit{, then}

\noindent (i) $S_{\circ }^{-1}=\{-s/(1+\rho s):s\in S\}\subseteq \mathbb{A}%
^{\varphi }$\textit{;}

\noindent (ii) $S\circ S=\{s+t\eta (s):s,t\in S\}\subseteq \mathbb{A}%
^{\varphi }$\textit{.}

\bigskip

\noindent \textit{Proof.}\textbf{\ }(i) As $S_{\circ }^{-1}$ is Baire and
non-meagre, Prop. 7 applies and $S_{\circ }^{-1}\subseteq \mathbb{A}_{\text{u%
}}\subseteq \mathbb{A}^{\varphi }$.

\noindent (ii) By Th. PJ, $S\circ S$ is isomorphic either to $S+S$ (for $%
\rho =0)$ or to $\eta _{\rho }(S)\eta _{\rho }(S)$ (for $\rho >0)$ and so is
Baire and non-meagre, by the Steinhaus Sum Theorem ([BinO9]); again Prop. 7
applies and $S\circ S\subseteq \mathbb{A}_{\text{u}}\subseteq \mathbb{A}%
^{\varphi }$. $\square $

\bigskip

\noindent \textit{Proof of Theorem 5. }Suffice it to assume\textit{\ }$\rho
=\rho _{\varphi }>0.$ Replacing $S$ by $S\cup (S_{\circ }^{-1})$ if
necessary, we may assume by Cor. 4 that $S$ is symmetric ($S=S_{\circ
}^{-1}) $, and w.l.o.g. $0=1_{\mathbb{G}}\in S,$ by Prop. 7.

Applying Cor. 4(ii) inductively, we deduce that%
\[
S^{\ast }:=\bigcup\nolimits_{n\in \mathbb{N}}(n)\circ S\subseteq \mathbb{A}%
^{\varphi },
\]%
where $(n)\circ S$ denotes $S\circ _{\eta }...\circ _{\eta }S$ to $n$ terms.
So $S^{\ast }$ is symmetric, and a semi-group: if $s\in (n)\circ S$ and $%
s^{\prime }\in (m)\circ S,$ then $s\circ s^{\prime }\in (n+m)\circ
S\subseteq S^{\ast }.$ So $\mathbb{A}^{\varphi }$ contains $S^{\ast }.$ As $%
0\in S^{\ast }$ (as above), $S^{\ast }$ is a subgroup (being symmetric,
since $\circ $ is commutative); hence $S^{\ast }$ is all of $\mathbb{G}%
_{+}^{\rho }.$ So $S^{\ast }=\mathbb{G}_{+}^{\rho }\mathbb{=A}_{\text{u}}=%
\mathbb{A}^{\varphi }.$ By Prop. 6, $\bar{K}(t)=K(\eta ^{-1}(e^{t}))$ is
additive on $\mathbb{R}$; indeed, by Prop. 6 with $\eta (u)=e^{x}$ and $\eta
(v)=e^{y}$%
\[
\bar{K}(x+y):=K(\eta ^{-1}(e^{x+y}))=K(u\circ _{\eta }v)=K(u)+K(v)=\bar{K}%
(x)+\bar{K}(y).
\]%
By Prop. 7 convergence is uniform near $u=0,$ so that $\bar{K}(t)$ is
bounded in a neighbourhood of $0,$ and, being additive, is linear; see e.g.
BGT 1.3, [Kuc], [BinO9,11]. So for some $c\in \mathbb{R}$ :
\[
c\log \eta (u)=c\log (1+\rho u)=cx=\bar{K}(x)=\bar{K}(\eta
^{-1}(e^{x}))=K(u)\quad (u>\rho ^{\ast }).\quad \square
\]

\section{Quantifier weakening}

Here we again drop the assumption that $\mathbb{A}^{\varphi }$ is co-meagre;
instead we will impose a density assumption, and employ a subadditivity
argument developed in [BinO12]. To motivate this, we recall the following
decomposition theorem of a function, with a one-sided finiteness condition,
into two parts, one decreasing, one with suitable limiting behaviour.

\bigskip

\noindent \textbf{Theorem BG2} ([BinG2, Th. 7]). \textit{The following are
equivalent:}

\noindent (i) \textit{The function }$U$\textit{\ has the decomposition}%
\[
U(x)=V(x)+W(x),
\]%
\textit{where }$V$\textit{\ has linear limiting moving average }$K_{V}$%
\textit{\ as in \S 1, and }$W(x)$\textit{\ is non-increasing;}

\noindent (ii) \textit{the following limit is finite:}%
\[
\lim_{\delta \downarrow 0}\limsup_{x\rightarrow \infty}\sup \left\{ \frac{%
U(x\circ _{\varphi }t)-U(x)}{\delta \varphi (x)}:t\in I_{\delta
}^{+}(0)\right\} <\infty .
\]

\noindent \textbf{Definitions. }For $\varphi \in SE$ and $\rho =\rho
_{\varphi },$ put%
\begin{eqnarray*}
H^{\dagger }(t) &:=&\lim_{\delta \downarrow 0}\limsup_{x\rightarrow \infty
}\sup \left\{ h(x\circ _{\varphi }s)-h(x):s\in I_{\delta }^{+}(t)\right\}
\qquad (t>\rho ^{\ast }), \\
\mathbb{A}_{\text{u}}^{\dagger } &:=&\{t>\rho ^{\ast }:H^{\dagger
}(t)<\infty \}.
\end{eqnarray*}%
So $\mathbb{A}_{\text{u}}\subseteq \mathbb{A}_{\text{u}}^{\dagger }$, as $%
H^{\dagger }(t)=K(t)$ on $\mathbb{A}_{\text{u}}$. In Theorem 6 below we
apply the techniques of [BinO11,12]; a first step for this is the following.
Here it is again convenient to rely on Prop. 5.

$\bigskip $

\noindent \textbf{Proposition 8. }\textit{For }$\varphi \in SE$\textit{\ and
}$\eta =\eta ^{\varphi },$\textit{\ }$H^{\dagger }\mathit{\ }$\textit{is
subadditive on }$\mathbb{A}_{\text{u}}^{\dagger }$ \textit{over non-inverse
pairs of elements }$s,t$\textit{:}
\[
H^{\dagger }(s\circ _{\eta }t)\leq H^{\dagger }(s)+H^{\dagger }(t)\mathit{%
\quad }(s,t\in \mathbb{A}_{\text{u}}^{\dagger },s\circ _{\eta }t\neq 0).
\]%
\textit{If }$\rho _{\varphi }^{\ast }$ \textit{is an accumulation point of }$%
\mathbb{A}_{\text{u}}^{\dagger }$\textit{, then either }$\lim
\inf_{s\downarrow \rho ^{\ast }}H^{\dagger }(s)$ \textit{is infinite, or }$%
H^{\dagger }\geq 0$\textit{\ on} $\mathbb{A}_{\text{u}}^{\dagger }.$

\bigskip

\noindent \textit{Proof.} For $c=\{c_{n}\}\rightarrow 1$ and $x=\{x_{n}\}$
divergent, put%
\[
H(t;x,c):=\lim \sup h(x_{n}\circ _{\varphi }c_{n}t)-h(x_{n})\mathit{\quad }%
(t\neq 0).
\]%
As in Prop. 6, for a given $c_{n}\rightarrow 1$ and divergent $x_{n},$ take $%
y_{n}:=x_{n}\circ _{\varphi }c_{n}s,$ $d_{n}:=c_{n}\eta (s)\varphi
(x_{n})/\varphi (y_{n})\rightarrow \eta (s)\eta (s)^{-1}=1.$ Now%
\[
x_{n}\circ _{\varphi }c_{n}(s+t\eta (s))=x_{n}+c_{n}(s+t\eta (s))\varphi
(x_{n})=y_{n}+d_{n}t\varphi (y_{n}),
\]%
so%
\[
h(x_{n}\circ _{\varphi }c_{n}(s+t\eta (s)))-h(x_{n})=h(y_{n}\circ _{\varphi
}d_{n}t)-h(y_{n})+h(x_{n}\circ _{\varphi }c_{n}s)-h(x_{n}),
\]%
whence
\[
H(s+t\eta (s);c,x)\leq H(t;d,y)+H(s;c,x).
\]%
So $H(s\circ _{\eta }t;c,x)\leq H^{\dagger }(t)+H^{\dagger }(s),$ as $%
H(t;d,y)\leq H^{\dagger }(t)$ and $H(s;c,x)\leq H^{\dagger }(s).$ Now we may
take suprema, since Prop. 5 applies provided $s\circ _{\eta }t\neq 0$.

For the final assertion, let $t\in \mathbb{A}_{\text{u}}^{\dagger }$ and
assume that $\lim \inf_{s\downarrow \rho ^{\ast }}H^{\dagger }(s)$ is
finite. Let $\varepsilon >0.$ Since $s\circ _{\eta }t\downarrow \rho ^{\ast
},$ for small enough $s$ with $\rho ^{\ast }<s<t_{\circ }^{-1}$
\[
-H^{\dagger }(t)+\lim \inf_{s\downarrow \rho ^{\ast }}H^{\dagger
}(s)-\varepsilon \leq H^{\dagger }(s\circ _{\eta }t)-H^{\dagger }(t)\leq
H^{\dagger }(s).
\]%
So%
\[
-H^{\dagger }(t)+\lim \inf_{s\downarrow \rho ^{\ast }}H^{\dagger
}(s)-\varepsilon \leq \lim \inf_{s\downarrow \rho ^{\ast }}H^{\dagger
}(s):\quad -\varepsilon \leq H^{\dagger }(t).
\]%
So $0\leq H^{\dagger }(t),$ as $\varepsilon >0$ was arbitrary$.$ $\square $

\bigskip

Our next result clarifies the role of the Heiberg-Seneta condition, for
which see BGT \S 3.2.1 and [BinO12]. It is here that we again use $(SE_{%
\mathbb{A}})$ on the set $\mathbb{A=G}_{+}^{\rho }=\{t:\eta _{\rho }(t)>0\}$
with $\rho =\rho _{\varphi }.$

\bigskip

\noindent \textbf{Proposition 9. }\textit{For }$\varphi \in SE,$ \textit{the
following are equivalent:}

\noindent (i) $0\in \mathbb{A}_{\text{u}}$ \textit{(i.e. }$\mathbb{A}_{\text{%
u}}\neq \emptyset $\textit{\ and so a subgroup);}

\noindent (ii) $\lim_{x\rightarrow \infty }[h(x+u\varphi (x))-h(x)]=0$
\textit{uniformly near} $u=0$\textit{;}

\noindent (iii) $H^{\dagger }(t)$\textit{\ satisfies the two-sided
Heiberg-Seneta condition:}%
\begin{equation}
\limsup_{u\rightarrow 0}H^{\dagger }(u)\leq 0.
\tag{$HS_{\pm
}(H^{\dagger
})$}
\end{equation}

\noindent \textit{Proof.} It is immediate that (i) and (ii) are equivalent.
We will show that (ii) and (iii) are equivalent. First assume the
Heiberg-Seneta condition. Take $\varepsilon >0,$ $x_{n}$ divergent, and $%
z_{n}$ null (i.e. $z_{n}\rightarrow 0$). By $HS_{\pm }(H^{\dagger })$, there
is $\delta _{\varepsilon }>0$ such that%
\[
H^{\dagger }(t)<\varepsilon \quad (0<|t|<\delta _{\varepsilon }).
\]%
So for each $t$ with $0<|t|<\delta _{\varepsilon }$ there are $\delta (t)>0$
and $X_{t}$ such that
\[
h(x\circ _{\varphi }s)-h(x)<\varepsilon \quad \quad (x>X_{t},s\in I_{\delta
(t)}^{+}(\pm t)).
\]%
By compactness, there are: $\delta >0,$ a finite set $F$ of points $t$ with $%
\delta _{\varepsilon }/3\leq t\leq 2\delta _{\varepsilon }/3,$ and $X$ such
that
\[
h(x\circ _{\varphi }s)-h(x)<\varepsilon \quad \quad (x>X,s\in I_{\delta
}^{+}(\pm t),t\in F),
\]%
and, further, $\{I_{\delta }^{+}(\pm t):t\in F\}$ covers $[-2\delta
_{\varepsilon }/3,-\delta _{\varepsilon }/3]\cup \lbrack \delta
_{\varepsilon }/3,2\delta _{\varepsilon }/3].$ By assumption, $\eta
_{x(n)}(s)\rightarrow 1$ uniformly as $s\rightarrow 0,$ so we may fix $%
t,t^{\prime }\in F$ and $s>0$ such that:\newline
\noindent (i) $s\in (t,t+\delta ),$\newline
\noindent (ii) w.l.o.g., for all $n,$ $s_{n}=s+z_{n}\in I_{\delta }^{+}(t),$%
\newline
\noindent (iii) w.l.o.g., for all $n,$ $-s/\eta _{x(n)}(s_{n})\in I_{\delta
}^{+}(-t^{\prime }).$\newline
\noindent Take $y_{n}=x_{n}\circ _{\varphi }s_{n}$; then for $x_{n},y_{n}>X,$
as in the proof of Theorem 3,%
\begin{eqnarray*}
&&h(x_{n}+z_{n}\varphi (x_{n}))-h(x_{n}) \\
&=&h(x_{n}+s_{n}\varphi (x_{n})-s\varphi (x_{n}))-h(x_{n}) \\
&=&h(y_{n}-s/\eta _{x(n)}(s_{n})\varphi (y_{n}))-h(y_{n})+h(x_{n}\circ
_{\varphi }s_{n})-h(x_{n})\leq 2\varepsilon .
\end{eqnarray*}%
In summary: for any divergent $x_{n}$ and null $z_{n}$%
\begin{equation}
h(x_{n}+z_{n}\varphi (x_{n}))-h(x_{n})<2\varepsilon \qquad \text{for all
large }n.  \label{2eps}
\end{equation}%
Towards a similar lower bound, suppose that for some divergent $y_{n}$ and
null $z_{n}^{\prime }$%
\[
h(y_{n}\circ _{\varphi }z_{n}^{\prime })-h(y_{n})\leq -2\varepsilon \qquad
\text{for all }n.
\]%
Take $x_{n}:=y_{n}\circ _{\varphi }z_{n}^{\prime },$ which is divergent;
then $y_{n}=x_{n}\circ _{\varphi }z_{n}$ for $z_{n}:=-z_{n}^{\prime }\varphi
(y_{n})/\varphi (x_{n})$ which is null$,$ since $\varphi (y_{n})/\varphi
(y_{n}\circ _{\varphi }z_{n}^{\prime })\rightarrow 1$ (by locally uniform
convergence of $\eta _{y(n)}$ near $0).$ So for all $n$%
\[
h(x_{n})-h(x_{n}+z_{n}\varphi (x_{n}))\leq -2\varepsilon :\qquad
h(x_{n}+z_{n}\varphi (x_{n}))-h(x_{n})\geq 2\varepsilon ,
\]%
a contradiction to (\ref{2eps}) for $n$ large enough.

So the Heiberg-Seneta condition yields%
\[
\lim [h(x+u\varphi (x))-h(x)]=0\text{ uniformly near }u=0,
\]%
i.e. (ii) holds.

Conversely, assuming (ii), for given $\varepsilon >0$ there are $X>0$ and $%
d>0$ so that for $x>X$ and $|u|<d,h(x+u\varphi (x))-h(x)<\varepsilon .$ So
for $x>X$
\[
\sup \{h(x+u\varphi (x))-h(x):|u|<d\}\leq \varepsilon .
\]%
Fixing $t\in (-d,d),$ choose $\delta >0$ so small that $I_{\delta
}(t)\subseteq (-d,d);$ then
\[
H_{\delta }(t):=\limsup_{x\rightarrow \infty }\sup \{h(x+u\varphi
(x))-h(x):u\in I_{\delta }^{+}(t)\}\leq \varepsilon .
\]%
But $H_{\delta }(t)$ is decreasing with $\delta ;$ so $H^{\dagger
}(t)=\lim_{\delta \downarrow 0}H_{\delta }(t)\leq \varepsilon $ for $t\in
(-d,d),$ i.e. $\lim \sup_{u\rightarrow 0}H^{\dagger }(u)\leq 0.$ $\square $

\bigskip

The final result of this section is the Beurling version of a theorem proved
in the Karamata framework of [BinO12]. However, uniformity plays no role
there, whereas here it is critical. The result shows that weakening the
quantifier in the definition of additivity to range only over a dense
subgroup, determined by locally uniform limits, yields `linearity' of $%
H^{\dagger }$. The $K$ in Th. 6 below is as in $(K)$ of \S 5, cf. Prop. 6.

\bigskip

\noindent \textbf{Theorem 6 (Quantifier Weakening from Uniformity).} \textit{%
For }$\rho >0$\textit{, if }$\mathbb{A}_{\text{u}}$\textit{\ is dense in }$%
\mathbb{G}_{+}^{\rho }$\textit{\ and } $H^{\dagger }(t)=K(t)$\textit{\ on }$%
\mathbb{A}_{\text{u}}$\textit{\ -- i.e. }$H^{\dagger }:(\mathbb{A}_{\text{u}%
},\circ _{\rho })\rightarrow (\mathbb{R},+)$ \textit{is a homomorphism --
then }$\mathbb{A}_{\text{u}}=\mathbb{G}_{+}^{\rho }$ \textit{and for some }$%
c\in \mathbb{R}$\textit{: }%
\[
H^{\dagger }(t)=c\log (1+\rho t)\quad \text{(}t>\rho ^{\ast }).
\]

\noindent \textit{Proof.}\textbf{\ }We check that Theorem 1 of [BinO12]
applies respectively to $\bar{H}(t):=H^{\dagger }(\eta ^{-1}(e^{t}))$ and $%
\bar{K}(t):=K(\eta ^{-1}(e^{t}))$ in place of $H$ and $K$ there, and with $%
\mathbb{A}:=\eta ^{-1}(\exp [\mathbb{A}_{\text{u}}]),$ which is dense in $%
\mathbb{R}$, since $\eta $ is an isomorphism taking $(\mathbb{G}_{+}^{\rho
},\circ _{\rho })$ to $(\mathbb{R}_{+},\cdot )$ (by Theorem PJ). Indeed, as
in Theorem 5 $\bar{K}$ is additive on $\mathbb{R}$ (by Prop. 6), and
likewise, by Prop. 8, $\bar{H}$ is subadditive. As $e^{0}=1=\eta (0),$ $\bar{%
H}$ satisfies the Heiberg-Seneta condition, by Prop. 9. Finally, since $%
H^{\dagger }(t)=K(t)$ on $\mathbb{A}_{\text{u}}$, $\bar{H}(t)=\bar{K}(t)$ on
$\mathbb{A}$. So $\bar{K}$ is linear by [BinO12, Th. 1], and the conclusion
follows once again as in Theorem 5. $\square $

\bigskip

\textbf{Remark.} As $\log [1+\rho (u\circ _{\rho }v)]=\log [(1+\rho
u)(1+\rho v)],$ the function $c\log (1+\rho t)$ is `subadditive' in the
sense of Prop. 8 (indeed, perhaps `additive').

\section{Representation}

We begin by identifying the limiting moving average $K_{F}$ of \S 1. Below $%
\varphi ,$ being increasing, is Baire.

\bigskip

\noindent \textbf{Lemma 3.} \textit{If }$\varphi \in SE$\textit{\ is
increasing and the following limit exists for }$F:\mathbb{R}\rightarrow
\mathbb{R}$\textit{:}%
\[
K_{F}(u):=\lim \frac{F(x\circ _{\varphi }u)-F(x)}{\varphi (x)},\qquad
(u>\rho _{\varphi }^{\ast })
\]%
\textit{-- then }$K_{F}$ \textit{as above satisfies}%
\[
K_{F}(u\circ _{\eta }v)=K_{F}(u)+K_{F}(v)\eta (u)\text{ for }\eta =\eta
^{\varphi };
\]%
\textit{if }$F$ \textit{is Baire/measurable, then }$K_{F}$\textit{\ and }$%
\eta =\eta ^{\varphi }$\textit{\ are of the form }%
\[
K_{F}(u)=c_{F}u,\text{ }\eta (u)=1+\rho u.
\]

\noindent \textit{Proof.}\textbf{\ }Write $y=x+u\varphi (x);$ then $\varphi
(y)/\varphi (x)\rightarrow \eta (u)$. Now
\begin{eqnarray*}
\frac{F(x\circ _{\varphi }[u+v])-F(x)}{\varphi (x)} &=&\frac{F(y\circ
_{\varphi }[v\varphi (x)/\varphi (y)])-F(y)}{\varphi (y)}\frac{\varphi (y)}{%
\varphi (x)} \\
&&+\frac{F(x\circ _{\varphi }u)-F(x)}{\varphi (x)}.
\end{eqnarray*}%
Write $w:=v/\eta (u);$ then, taking limits above, gives%
\[
K_{F}(u+w\eta (u))=K_{F}(w)\eta (u)+K_{F}(u).
\]%
Assuming $F$ is Baire/measurable, $K_{F}(t)=\lim_{n\rightarrow \infty
}[(F(n\circ _{\varphi }u)-F(x))/\varphi (n)]$ is Baire/measurable (as in
Prop. 7). By [BinO11, Th.9,10] $K_{F}(x)=c_{F}H_{0}(x),$ where $H_{0}(x):=x.$
So $K_{F}(u)=c_{F}u,$ for some $c_{F}.\square $

\bigskip

The result above formally extends to (i) the Beurling framework, and (ii) to
the class $SE$ the notion of $\Pi _{g}$\textit{-}class, due to Bojani\'{c}%
-Karamata/de Haan, for which see BGT Ch. 3, since just as there
\begin{equation}
\text{(i) }\frac{F(x\circ _{\varphi }u)-F(x)}{\varphi (x)}\sim
c_{F}H_{0}(u):\qquad \text{(ii) }\frac{F(x\circ _{\varphi }u)-F(x)}{u\varphi
(x)}\rightarrow c_{F}.  \tag{$\Pi _{\varphi }$}
\end{equation}%
\noindent \textbf{Definition. }Say that $F$ is of \textit{Beurling} $\Pi
_{\varphi }$\textit{-class} \textit{with }$\varphi $\textit{-index }$c=c_{F}$%
\textit{\ }(cf. BGT Ch. 3) if the convergence in $(\Pi _{\varphi }$(ii)$)$
is locally uniform in $u.$

\bigskip

This should be compared with Theorem BG2 in \S 7. We now use a Goldie-type
argument (see [BinO11]) to establish the representation below for the class $%
\Pi _{\varphi }$.

\bigskip

\noindent \textbf{Theorem 7 (Representation for }Beurling $\Pi _{\varphi }$%
-class with $\varphi $-index $c$). \textit{For }$F$\textit{\
Baire/measurable, }$F$ \textit{is of additive Beurling} $\Pi _{\varphi }$%
\textit{-class} \textit{with }$\varphi $\textit{-index }$c$ \textit{and }$%
\varphi \in SE$ \textit{iff}%
\[
F(x)=b+cx+\int_{1}^{x}e(t)dt,\text{ }b\in \mathbb{R}\text{ and }e\rightarrow
0.
\]

\noindent \textit{Proof.} As above, by the $\lambda $-UCT of [Ost3, Th. 1],\
there exists $X$ such that for all $x\geq X$ and all $u$ with $|u|\leq 1$%
\[
\frac{F(x\circ _{\varphi }u)-F(x)}{u\varphi (x)}=c+\varepsilon (x;u),
\]%
with%
\[
\varepsilon (x;u)\rightarrow 0\text{ uniformly for }|u|\leq 1\text{ as }%
x\rightarrow \infty .
\]%
Put%
\[
e(x):=\sup \{\varepsilon (x,u):|u|\leq 1\};
\]%
then $e(x)\rightarrow 0$ as $x\rightarrow \infty .$

Using a Beck sequence $x_{n+1}=x_{n}\circ _{\varphi }u$ ([BinO11, \S 3]; cf.
Bloom [Blo], BGT\ Lemma 2.11.2) starting at $x_{0}=X$ and ending at $%
x_{m}=x(u)\leq x$ with $x<x(u)\circ _{\varphi }u=x_{m+1}$ yields
\begin{eqnarray*}
F(x(u))-F(X)
&=&\sum\nolimits_{n=0}^{m-1}F(x_{n+1})-F(x_{n})=\sum\nolimits_{n=0}^{m-1}(c+%
\varepsilon (x_{n};u))u\varphi (x_{n}) \\
&=&\sum\nolimits_{n=0}^{m-1}(c+\varepsilon (x_{n};u))(x_{n}+u\varphi
(x_{n})-x_{n}) \\
&=&c\sum\nolimits_{n=0}^{m-1}(x_{n+1}-x_{n})+\sum\nolimits_{n=0}^{m-1}%
\varepsilon (x_{n};u)(x_{n+1}-x_{n}) \\
&=&c(x(u)-X)+\sum\nolimits_{n=0}^{m-1}\varepsilon (x_{n};u)(x_{n+1}-x_{n}).
\end{eqnarray*}%
Since $F$ is Baire/measurable, we may restrict attention to points $x$ where
$F$ is continuous. For $x$ fixed, note that $u\varphi (x_{n})\leq u\varphi
(x)\rightarrow 0$ as $u\rightarrow 0,$ so $x(u)\rightarrow x;$ taking limsup
as $u\rightarrow 0$,%
\[
F(x)=F(X)+c(x-X)+\int_{X}^{x}e(t)dt,
\]%
with $e(x)\rightarrow 0,$ as above. So on differencing,%
\[
\frac{F(x+u\varphi (x))-F(x)}{u\varphi (x)}=c+\frac{1}{u\varphi (x)}%
\int_{x}^{x+u\varphi (x)}e(t)dt\rightarrow c.
\]%
So $F$ is Beurling $\Pi _{\varphi }$-class with $\varphi $-index $c$\ iff it
has the representation stated\textit{. }$\square $

\bigskip

We note also a generalization of Prop. 8 and Lemma 2, for which we need
notation (similar to that in \S 7) analogous to the Karamata $\Omega $ of
BGT\ \S 3.0 (cf. BGT Th. 3.3.2/3).

\bigskip

\noindent \textbf{Definitions. }For $\varphi \in SE$ and $\rho =\rho
_{\varphi },$ put%
\begin{eqnarray*}
\Omega _{h}^{\dagger }(t)&:=&\lim_{\delta \downarrow 0}\limsup_{x\rightarrow
\infty }\sup \left\{ (h(x\circ _{\varphi }s)-h(x))/\varphi (x):s\in
I_{\delta }^{+}(t)\right\} , \\
\mathbb{A}_{\Omega }^{\dagger }&:= &\{t> \rho ^{\ast }:\Omega _{h}^{\dagger
}(t)<\infty \}.
\end{eqnarray*}

\noindent \textbf{Proposition 8}$^{\prime }$\textbf{. }\textit{For }$\varphi
\in SE$\textit{\ and }$\eta =\eta ^{\varphi }$\textit{\ }$,$\textit{\ }$%
\Omega _{h}^{\dagger }\mathit{\ }$\textit{is }$\eta $\textit{-subadditive on
}$\mathbb{A}_{\Omega }^{\dagger }$\textit{:}
\[
\Omega _{h}^{\dagger }(s\circ _{\eta }t)\leq \Omega _{h}^{\dagger }(t)\eta
(s)+\Omega _{h}^{\dagger }(s)\quad (s,t\in \mathbb{A}_{\Omega }^{\dagger
},s\in (\rho _{\varphi }^{\ast },+\infty )).
\]

\noindent \textit{Proof.} For $c=\{c_{n}\}\rightarrow 1$ and $x=\{x_{n}\}$
divergent, put%
\[
\Omega _{h}^{\dagger }(t;x,c):=\lim \sup [h(x_{n}\circ _{\varphi
}c_{n}t)-h(x_{n})]/\varphi (x_{n}).
\]%
As in Prop. 6, for a given $c_{n}\rightarrow 1$ and divergent $x_{n},$ take%
\[
y_{n}=x_{n}\circ _{\varphi }c_{n}s,\text{ }d_{n}:=1/\eta
_{x(n)}(s)\rightarrow \eta (s)^{-1}>0.
\]%
Since by Prop. 2(iii)%
\begin{eqnarray*}
&&[h(x_{n}\circ _{\varphi }c_{n}(s\circ _{\eta }t))-h(x_{n})]/\varphi (x_{n})
\\
&=&[h(y_{n}\circ _{\varphi }d_{n}t\eta (s))-h(y_{n})]/\varphi (y_{n})\cdot
\eta _{x(n)}(c_{n}s)+[h(x_{n}\circ _{\varphi }c_{n}s)-h(x_{n})]/\varphi
(x_{n}),
\end{eqnarray*}%
\[
\Omega _{h}^{\dagger }(s+t\eta (s);c,x)\leq \Omega _{h}^{\dagger
}(t;d,y)\eta (s)+\Omega _{h}^{\dagger }(s;c,y).
\]%
Now take suprema. $\square $

\bigskip

We note an extension of [BinG3, Th. 1] -- cf. the more recent [Bin].

\bigskip

\noindent \textbf{Theorem BG 3.} \textit{If }$\varphi \in SE$\textit{\ and }$%
\varphi \uparrow \infty ,$\textit{\ then }$U$\textit{\ has a limiting moving
average }$K_{U}(x)=cx$\textit{\ iff }%
\[
\frac{1}{\lambda (x)}\int_{0}^{x}U(y)d\lambda (y)\rightarrow c,
\]%
\textit{where }$\lambda (x):=\varphi (x)\exp \tau _{\varphi }(x).$

\bigskip

\noindent \textbf{Corollary 5. }\textit{For} $\varphi \in SE$\textit{\ and }$%
\varphi \uparrow \infty ,$\textit{\ and with }$\lambda $\textit{\ as
previously, if }$F$ \textit{is of additive Beurling} $\Pi _{\varphi }$%
\textit{-class} \textit{with }$\varphi $\textit{-index }$c,$ \textit{then}%
\[
\frac{1}{\lambda (x)}\int_{0}^{x}F(y)d\lambda (y)\rightarrow c.
\]

\section{Divided difference and double sweep}

The concern of previous sections was the asymptotics of differences: $\Delta
_{t}^{\varphi }h$ in the Beurling theory, and exceptionally in \S 8 moving
averages $\Delta _{t}^{\varphi }h/\varphi $ in the Beurling version of the
Bojani\'{c}-Karamata/de Haan theory. Introducing an appropriate general
denominator $\psi $ carries the same advantage as in BGT (e.g. 3.13.1) of
`double sweep': capturing the former theory via $\psi \equiv 1$ and the
latter via $\psi \equiv \varphi ,$ embracing both through a common
generalization -- see Prop. 8$^{\prime }$ above for a first hint of such
possibilities. The work of this section is mostly to identify how earlier
results generalize, much of it focussed on \S 3, to which we refer for
group-theoretic notation; in particular $\mathbb{G}$ denotes the relevant
\textit{Popa group}, i.e. $\mathbb{G}^{\rho }$ for $\rho =\rho _{\varphi }$
for the appropriate $\varphi ,$ with $\mathbb{G}_{+}:=\{t:\eta _{\rho
}(t)>0\}$ its positive half-line.

Let $\varphi \in SE;$ fix a $\varphi $-regularly varying $\psi >0$ with $%
\varphi $-index $\gamma $ and limit function $g,$ i.e.
\begin{equation}
\psi (x+t\varphi (x))/\psi (x)\rightarrow g(t)\text{ loc. uniformly in }t%
\text{ \qquad (}t>\rho ^{\ast }\text{),}  \tag{$G$}
\end{equation}%
and, since $g(t)$ is a homomorphism (see Prop. 10 below), it is either $%
e^{\gamma t}$ ($\rho =0),$ or else $\eta _{\rho }(t)^{\gamma }$ (see
[Ost4]). Recalling the notation $\Delta _{t}^{\varphi }h(x)$ from \S 1, we
also write $\Delta _{t}^{\varphi }h/\psi (x)$ to mean $(\Delta _{t}^{\varphi
}h(x))/\psi (x).$ We are concerned below with
\begin{equation}
H^{\ast }(t):=\limsup_{x}[\Delta _{t}^{\varphi }h/\psi ],  \tag{$H^{\ast }$}
\end{equation}%
whenever this exists, and with the nature of the convergence. To specify
whenever a case below of convergence arises, we write%
\begin{eqnarray*}
H_{+}^{\ast }(t){}:= &&\lim_{\delta \downarrow 0}\limsup_{x\rightarrow
\infty }\sup \{\Delta _{s}^{\varphi }h/\psi :s\in I_{\delta }^{+}(t)\}\text{,%
} \\
H_{-}^{\ast }(t){}:= &&\lim_{\delta \downarrow 0}\limsup_{x\rightarrow
\infty }\sup \{\Delta _{s}^{\varphi }h/\psi :s\in I_{\delta }^{-}(t)\}\text{,%
} \\
H_{\pm }^{\ast }(t){}:= &&\lim_{\delta \downarrow 0}\limsup_{x\rightarrow
\infty }\sup \{\Delta _{s}^{\varphi }h/\psi :s\in I_{\delta }(t)\}\text{.}
\end{eqnarray*}

We begin with an extension of Lemma 2, for which we recall the notation $%
\mathbb{A}_{\Omega }^{\dagger }$ of \S 8. The proofs are almost identical --
so are omitted.

\bigskip

\noindent \textbf{Lemma 2}$^{\dagger }$\textbf{.}\textit{\ }(i) \textit{If }$%
\varphi \in SE$ \textit{and }$(G)$ \textit{holds -- then:}

\noindent (a)\textit{\ if the convergence in }$(H^{\ast })$ \textit{is
uniform (resp. right-sidedly uniform) near }$t=0$\textit{, then it is
uniform (resp. right-sidedly uniform) everywhere in }$\mathbb{A}_{\Omega
}^{\dagger }$ \textit{and for }$u\in \mathbb{A}_{\Omega }^{\dagger }\cap
\mathbb{R}_{+}$\textit{\ }%
\[
H_{+}^{\ast }(u)\leq H^{\ast }(u)g(u)+H_{+}^{\ast }(0);
\]

\noindent (b) \textit{if the convergence in }$(H^{\ast })$\textit{\ is
uniform near }$t=u\in \mathbb{A}_{\Omega }^{\dagger }\cap (\mathbb{A}%
_{\Omega }^{\dagger })_{\circ }^{-1}\cap \mathbb{R}_{+},$\textit{\ then it
is uniform near }$t=0:$%
\[
H_{\pm }^{\ast }(0)\leq H_{\pm }^{\ast }(u)g(u_{\circ }^{-1})+H^{\ast
}(u_{\circ }^{-1});
\]

(ii) \textit{if }$\rho =0$ \textit{and }$\varphi \in SN$ \textit{is
monotonic increasing and the convergence in }$(H^{\ast })$\textit{\ is
right-sidedly uniform near }$t=u\in \mathbb{A}_{\Omega }^{\dagger }\cap
\mathbb{R}_{+},$\textit{\ then it is right-sidedly uniform near }$t=0:$%
\[
H_{+}^{\ast }(0)\leq H_{+}^{\ast }(u)g(u_{\circ }^{-1})+H^{\ast }(u_{\circ
}^{-1}).
\]

For the next result, recall also the notation $\Omega _{h}^{\dagger }(t)$ of
\S 8.

\bigskip

\noindent \textbf{Proposition 10. }(i)\textit{\ With }$g$\textit{\ as in }$%
(G)$ \textit{above, }$g(u\circ _{\eta }v)=g(u)g(v),$ \textit{so that}%
\[
K_{h}(u\circ _{\eta }v\circ _{\eta }w)=K_{h}(u)g(v\circ _{\eta
}w)+K_{h}(v)g(w)+K_{h}(w),
\]%
\textit{\ and furthermore}
\[
H^{\ast }(s\circ _{\eta }t)\leq H^{\ast }(t)g(s)+H^{\ast }(s)\text{ \qquad }%
(s,t\in \mathbb{A}_{\Omega }^{\dagger }).
\]

\noindent (ii)\textit{\ If both of the following hold:}

(a) $H^{\ast }(t)>-\infty $\textit{\ for }$t$\textit{\ in a subset }$\Sigma $%
\textit{\ that is unbounded below;}

(b) \textit{the Heiberg-Seneta condition }$\Omega _{h}^{\dagger }(0+)\leq 0$%
\textit{\ holds}

\noindent \textit{--- then }$H^{\ast }$\textit{\ is finite on }$\mathbb{G}%
_{+}$ \textit{and }$H^{\ast }(0+)=0.$

\noindent \textit{\ Moreover, for }$\mathbb{A}_{\Omega }^{\dagger }$\textit{%
\ dense in }$\mathbb{G}_{+}$\textit{,}%
\[
H^{\ast }(u\circ _{\eta }v)=K(v)g(u)+H^{\ast }(u)\text{ \qquad }(u\in
\mathbb{G}_{+}^{\rho },v\in \mathbb{A}^{\varphi }).
\]

\bigskip

\noindent \textit{Proof.}\textbf{\ }(i) The first assertion follows by
writing $y=x\circ _{\varphi }u$ (as in Prop. 2(iii)) and taking limits in
the identity
\[
\psi (x\circ _{\varphi }(u\circ _{\eta }v))/\psi (x)=[\psi (y\circ _{\varphi
}v/\eta _{x}(u)))/\psi (y)]\psi (x\circ _{\varphi }u)/\psi (x).
\]%
The assertion is a restatement of the Cauchy exponential equation for $%
e^{\gamma x}$ when $\rho =0$ and for $\eta (x)^{\gamma }$ for $\rho >0,$ and
so implies the second. As for the third assertion, argue as in Prop. 8$%
^{\prime }$ above, but now with a new denominator $\psi (x_{n})$.

(ii) The first assertion is proved from (a) as in [BinO12, Prop. 6], and the
second from part (b) as in [BinO12, Prop. 8]; the latter uses part (i) and
the two facts that $g(u\circ v)=g(u)g(v)$ and $g(u)\geq 1$ for $u>0.$ The
second assertion is proved as in BGT Th. 3.2.5. $\square $

\bigskip

As a corollary, since $H^{\ast }$ is $g$-subadditive, we have the analogue
of Th. 1 of [BinO12].

\bigskip

\noindent \textbf{Theorem 8. }\textit{In the setting of Proposition 10, if }$%
\mathbb{A}_{\Omega }^{\dagger }$\textit{\ is dense, then }$\mathbb{A}%
_{\Omega +}^{\dagger }=\mathbb{G}_{+}$\textit{\ and for some }$c,\gamma
,\rho \in \mathbb{R}$:

\noindent either\textit{\ }(i)\textit{\ }$\rho =0$\textit{\ and }$H^{\ast
}(u)\equiv cH_{(-\gamma )}(u)=c(1-e^{-\gamma u})/\gamma $\textit{\ }$(u\geq
0),$

\noindent or (ii)\textit{\ }$\rho >0$ \textit{and }$H^{\ast }(u)\equiv
\lbrack (1+\rho u)^{\gamma +1}-1]/[\rho (1+\gamma )]$\textit{\ }$(u\geq 0).$

\bigskip

\noindent \textit{Proof.}\textbf{\ }As in Prop. 6 above, $(\mathbb{A}%
_{\Omega }^{\dagger },\circ )$ is a subgroup. Now use Prop. 10, Theorem PJ,
and Th. 3 of [BinO11]. $\square $

\section{Uniform Boundedness Theorem}

As above, let $h$ be Baire and $\varphi \in SE$ on $\mathbb{R}_{+}$ be
positive. Thus for all divergent $x_{n}$ (i.e. divergent to +$\infty ),$
\[
\varphi (x_{n}\circ _{\varphi }t)/\varphi (x_{n})\rightarrow \eta (t)\text{
for all }t\geq 0\text{ and }\varphi (x)=O(x).
\]%
So $y_{n}=x_{n}\circ _{\varphi }t=x_{n}(1+t\varphi (x_{n})/x_{n})$ is
divergent if $x_{n}$ is.

We work additively, and recall that for $t>\rho ^{\ast }$%
\[
H^{\ast }(t):=\limsup_{x\rightarrow \infty }h(x\circ _{\varphi }t)-h(x),%
\text{ }H_{\ast }(t):=\liminf_{x\rightarrow \infty }h(x\circ _{\varphi
}t)-h(x).
\]%
If $x_{n}\rightarrow \infty $ and $H^{\ast }(t)<\infty ,$ then for all large
enough $n$%
\[
h(x_{n}\circ _{\varphi }t)-h(x_{n})<n.
\]%
Likewise if $H_{\ast }(t)>-\infty ,$ then for all large enough $n$%
\[
h(x_{n})-h(x_{n}\circ _{\varphi }t)<n.
\]

In the theorem below we need to assume finiteness of both $H^{\ast }$ and $%
H_{\ast };$ we recall that in the Karamata case, substituting $y$ for $u+x,$
one has%
\[
h^{\ast }(u)=\lim \sup [h(u+x)-h(x)]=-\lim \inf [h(y-u)-h(y)]=-h_{\ast
}(-u).
\]%
This relationship is used implicitly in the standard development of the
Karamata theory -- see e.g. BGT, \S 2.1. Theorem 9 below extends [BinO8, Th.
8]. As the hypothesis is symmetric, the same proof yields the liminf case.

\bigskip

\noindent \textbf{Theorem 9 (Uniform Boundedness Theorem}; cf. [Ost2]).
\textit{For }$\varphi \in SE$ \textit{and }$\rho =\rho _{\varphi }$\textit{,
suppose that }$-\infty <H_{\ast }(t)\leq H^{\ast }(t)<\infty $\textit{\ for}
$t\in S$\textit{\ with }$S\subseteq \mathbb{G}_{+}^{\rho }$ \textit{a
non-meagre Baire set.} \textit{Then for compact }$K\subseteq S$%
\[
\limsup_{x\rightarrow \infty }\left( \sup_{u\in K}h(x\circ _{\varphi
}u)-h(x)\right) <\infty .\text{ }
\]

\bigskip

\noindent \textit{Proof.} By compactness of $K,$ it suffices to establish
uniform boundedness locally at any point $u>\rho ^{\ast }$. Suppose
otherwise, and that this is witnessed by some $x_{n}\rightarrow \infty $ and
$u_{n}\rightarrow u.$ Writing $u_{n}:=u+z_{n}$ with $z_{n}\rightarrow 0$ and
passing if necessary from $x_{n}$ to $\xi _{n}:=x_{n}\circ _{\varphi }u$
(and using the identity $h(x_{n}\circ _{\varphi }u_{n})-h(x_{n})=[h(\xi
_{n}\circ _{\varphi }z_{n}/\eta _{x}(u))$ $-h(\xi _{n})]+[h(x_{n}\circ
_{\varphi }u)-h(x_{n})],$ where the first bracket tends to 0) w.l.o.g. we
may assume $u=0,$ and%
\begin{equation}
h(x_{n}\circ _{\varphi }z_{n})-h(x_{n})>3n.  \label{big}
\end{equation}%
Put $y_{n}:=x_{n}\circ _{\varphi }z_{n}.$ As $\varphi \in SE,$%
\[
c_{n}:=\varphi (x_{n}\circ _{\varphi }z_{n})/\varphi (x_{n})\rightarrow 1.
\]%
Write $\gamma _{n}(s):=c_{n}s+z_{n}$ . Put
\begin{eqnarray*}
V_{n}:= &\{s\in &S:h(x_{n}\circ _{\varphi }s)-h(x_{n})<n\},\text{ }%
H_{k}^{+}:=\bigcap\nolimits_{n\geq k}V_{n}, \\
W_{n}:= &\{s\in &S:h(y_{n})-h(y_{n}\circ _{\varphi }s)<n\},\text{ }%
H_{k}^{-}:=\bigcap\nolimits_{n\geq k}W_{n}.
\end{eqnarray*}%
These are Baire sets, and since $-\infty <H_{\ast }(t)\leq H^{\ast
}(t)<\infty $ on $S,$%
\begin{equation}
S=\bigcup\nolimits_{k}H_{k}^{+}=\bigcup\nolimits_{k}H_{k}^{-}.  \label{cov}
\end{equation}%
The increasing sequence of sets $\{H_{k}^{+}\}$ covers $S.$ So for some $k$
the set $H_{k}^{+}$ is non-negligible. Then, by (\ref{cov}), for some $l$
the set%
\[
B:=H_{k}^{+}\cap H_{l}^{-}
\]%
is also non-negligible. Take $A:=H_{k}^{+};$ then $B\subseteq H_{l}^{-}$ and
$B\subseteq A$ with $A,B\ $non-negligible. Applying the Affine Two-sets
Lemma [BinO10, Lemma 2] to the maps $\gamma _{n}(s)=c_{n}s+z_{n}$ with $%
c=\lim_{n}c_{n}=1,$ there exist $b\in B$ and an infinite set\textit{\ }$%
\mathbb{M}$ with%
\[
\{c_{m}b+z_{m}:m\in \mathbb{M}\}\subseteq A=H_{k}^{+}.
\]%
That is, as $B\subseteq H_{l}^{-},$ there exist $t\in H_{l}^{-}$ and an
infinite $\mathbb{M}_{t}$ with%
\[
\{\gamma _{m}(t)=c_{m}t+z_{m}:m\in \mathbb{M}_{t}\}\subseteq H_{k}^{+}.
\]%
In particular, for this $t$ and $m\in \mathbb{M}_{t}$ with $m>k,l$,
\[
t\in W_{m}\hbox{ and }\gamma _{m}(t)\in V_{m}.
\]%
As $\gamma _{m}(t)\in V_{m},$%
\begin{equation}
h{(x_{m}\circ _{\varphi }\gamma }_{m}({t))}-h{(x_{m})<m}.  \label{*}
\end{equation}%
But $\gamma _{m}(t)=z_{m}+c_{m}t=z_{m}+t\varphi (y_{m})/\varphi (x_{m}),$ so
\[
x_{m}\circ _{\varphi }\gamma _{m}(t)=x_{m}+z_{m}\varphi (x_{m})+t\varphi
(y_{m})=y_{m}\circ _{\varphi }t.
\]%
So, by (\ref{*}),
\[
h{(y_{m}\circ _{\varphi }t)}-h{(x_{m})<m}.
\]%
But $t\in W_{m},$ so%
\[
h{(y_{m})-}h{(y_{m}\circ _{\varphi }t)<m}.
\]%
Combining these with (\ref{cov}) and (\ref{big}).
\[
3m<h(y_{m})-h(x_{m})\leq \{h(y_{m})-h(y_{m}\circ _{\varphi
}t)\}+\{h(y_{m}\circ _{\varphi }t)-h(x_{m})\}\leq 2m,
\]%
a contradiction. $\square $

\bigskip

As in the classical Karamata case, this result implies global bounds on $h$
-- see BGT Th. 2.0.1.

\bigskip

\noindent \textbf{Theorem 10.} \textit{In the setting of Theorem 9, for }$%
\varphi \in SE,$\textit{\ if the set }$S$\textit{\ on which }$H^{\ast }(t)$%
\textit{\ and }$H_{\ast }(t)$\textit{\ are finite contains a half-interval }$%
[a_{0},\infty )$\textit{\ with }$a_{0}>0$\textit{\ --- then there is a
constant }$C>0$ \textit{such that for all large enough }$x$\textit{\ and }$u$
\[
h(u\varphi (x)+x)-h(x)\leq C\log u.
\]%
The proof parallels the end of the proof in BGT\ of Th. 2.0.1, but with the
usual sequence of powers $a^{n}$ replaced by a Popa-style generalization
(cf. Prop. 2(v)):%
\[
a_{\varphi x}^{n+1}:=a_{\varphi x}^{n}\circ _{\varphi x}a=a_{\varphi
x}^{n}+a\eta _{x}^{\varphi }(a_{\varphi x}^{n})\text{ with }a_{\varphi
x}^{1}=a.
\]%
It relies on estimation results for $a_{\varphi x}^{m}$ that are uniform in $%
m$ (this only needs $\eta _{x}^{\varphi }\rightarrow \eta _{\rho }$
pointwise):\footnote{%
See the Appendixfor the proofs of Th. 10 and Prop. 11.}

\bigskip

\noindent \textbf{Proposition 11.} \textit{If }$\varphi \in SE$\textit{\
with }$\rho =\rho _{\varphi }>0,$ \textit{then for any }$a>1,$ $%
0<\varepsilon <1,$\textit{\ }\newline
\noindent (i)($a_{\varphi x}^{m}$-estimates under $\eta _{x}^{\varphi }$)%
\textit{\ for all} \textit{large enough }$x$:%
\[
(1-\varepsilon )\leq \eta _{x}^{\varphi }(a_{\varphi x}^{m})^{1/m}/\eta
_{\rho }(a)\leq (1+\varepsilon ),\qquad (m\in \mathbb{N});
\]%
\noindent (ii)($a_{\varphi x}^{m}$-estimates under $\eta _{\rho }$) \textit{%
for all} \textit{large enough }$x$:\textit{\ }%
\[
\frac{\eta _{\rho }(a(1-\varepsilon ))^{m}}{1-\varepsilon }-\frac{%
\varepsilon }{1-\varepsilon }\leq \eta _{\rho }(a_{\varphi x}^{m})\leq \frac{%
\eta _{\rho }(a(1+\varepsilon ))^{m}}{1+\varepsilon }+\frac{\varepsilon }{%
1+\varepsilon },\qquad (m\in \mathbb{N});
\]%
\noindent (iii) $a_{\varphi x}^{m}\rightarrow \infty ;$ \textit{and\newline
\noindent }(iv)\textit{\ there are }$C_{\pm }=C_{\pm }(\rho ,a,\varepsilon
)>0$\textit{\ such that for all} \textit{large enough }$x$ \textit{and }$u$:%
\[
a_{\varphi x}^{m}\leq u<a_{\varphi x}^{m+1}\Longrightarrow mC_{-}\leq \log
u\leq (m+1)C_{+}.
\]

\section{Character degradation from limsup}

We refer the reader to [BinO8] for a discussion, from the perspective of the
practising analyst (employing `naive' set theory), of the broader
set-theoretic context below. For convenience we repeat part of the
commentary in [BinO8]; for more detail see [BinO13]. As there so too here,
our interest in the complexities induced by the \textit{limsup} operation
points us in the direction of definability and descriptive set theory,
because of the question of whether certain specific sets, encountered in the
course of the analysis, have the Baire property. The answer depends on what
further axioms one admits. For us there are two alternatives yielding the
kind of decidability we seek: G\"{o}del's Axiom of Constructibility $V=L$,
as an appropriate \textit{strengthening} of the Axiom of Choice (AC) which
creates definable sets without the Baire property (without measurability),
or, at the opposite pole, the Axiom of Projective Determinacy, $PD$ (see
[MySw], or [Kec, 5.38.C]), an \textit{alternative} to $AC$ which guarantees
the Baire property in the kind of definable sets we encounter. Thus to
decide whether sets of the kind we encounter below have the Baire property,
or are measurable, the answer is: it depends on the axioms of set theory
that one adopts. It turns out that $AC$ may be usefully weakened to the
Axiom of Dependent Choice(s), $DC;$ for details see [BinO13].

To formulate our results we need the language of descriptive set theory, for
which see e.g. [JayR], [Kec], [Mos]. Within such an approach we will regard
a function as a set, namely its \textit{graph}; formulas written in naive
set-theoretic notation then need a certain amount of formalization -- for a
quick approach to such matters refer to [Dra, Ch. 1,2] or the very brief
discussion in [Kun, \S 1.2]. We need the beginning of the \textit{projective
hierarchy }in Euclidean space (see [Kec, S. 37.A]), in particular the
following classes:

the \textit{analytic} sets $\mathbf{\Sigma }_{1}^{1};$

their complements, the \textit{co-analytic} sets $\mathbf{\Pi }_{1}^{1};$

the common part of the previous two classes, the \textit{ambiguous} class $%
\mathbf{\Delta }_{1}^{1}:=\mathbf{\Sigma }_{1}^{1}\cap \mathbf{\Pi }%
_{1}^{1}, $ that is, by Souslin's Theorem ([JayR, p. 5], and [MaKe, p.407]
or [Kec, 14. C]) the \textit{Borel} sets;

the \textit{projections} (continuous images) of $\mathbf{\Pi }_{1}^{1}$
sets, forming the class $\mathbf{\Sigma }_{2}^{1};$

their \textit{complements}, forming the class $\mathbf{\Pi }_{2}^{1};$

the \textit{ambiguous} class $\mathbf{\Delta }_{2}^{1}:=\mathbf{\Sigma }%
_{2}^{1}\cap \mathbf{\Pi }_{2}^{1}$;

and then: $\mathbf{\Sigma }_{n+1}^{1},$ the projections of $\mathbf{\Pi }%
_{n}^{1};$ their complements $\mathbf{\Pi }_{n+1}^{1};$ and the ambiguous
class $\mathbf{\Delta }_{n+1}^{1}:=\mathbf{\Sigma }_{n+1}^{1}\cap \mathbf{%
\Pi }_{n+1}^{1}.$

Throughout we shall be concerned with the cases $n=1,2$ or $3.$

The notation reflects the fact that the canonical expression of the logical
structure of their definitions, namely with the quantifiers (ranging over
the reals, hence the superscript 1, as reals are type 1 objects -- integers
are of type 0) all at the front, is determined by a string of alternating
quantifiers starting with an existential or universal quantifier (resp. $%
\mathbf{\Sigma }$ or $\mathbf{\Pi })$. Here the subscript accounts for the
number of alternations.

Interest in the character of a function $H$ is motivated by an interest
within the theory of regular variation in the character of the level sets%
\[
H^{k}:=\{s:|H(s)|<k\}=\{s:(\exists t)[(s,t)\in H\text{ }\&\text{ }|t|<k]\},
\]%
for $k\in \mathbb{N}$ (where as above $H$ is identified with its graph). The
set $H^{k}$ is thus the projection of $H\cap (\mathbb{R}\times \lbrack
-k,k]) $ and hence is $\mathbf{\Sigma }_{n}^{1}$ if $H$ is $\mathbf{\Sigma }%
_{n}^{1},$ e.g. it is $\mathbf{\Sigma }_{1}^{1},$ i.e. analytic, if $H$ is
analytic (in particular, Borel). Also
\[
H^{k}=\{s:(\forall t)[(s,t)\in H\text{ }\Longrightarrow \text{ }|t|\leq
k]\}=\{s:(\forall t)[(s,t)\notin H\text{ or }|t|\leq k]\},
\]%
and so this is also $\mathbf{\Pi }_{n}^{1}$ if $H$ is $\mathbf{\Sigma }%
_{n}^{1}.$ Thus if $H$ is $\mathbf{\Sigma }_{n}^{1}$ then $H^{k}$ is $%
\mathbf{\Delta }_{n}^{1}.$ So if $\mathbf{\Delta }_{n}^{1}$ sets are Baire,
then for some $k$ the set $H^{k}$ is Baire non-null.

With this in mind, it suffices to consider upper limits; as before, we
prefer to work with the additive formulation. Consider the definition:%
\begin{equation}
H_{\varphi }^{\ast }(x):=\lim \sup_{t\rightarrow \infty }[h(t+x\varphi
(t))-h(t)].  \tag{$\ast \ast $}
\end{equation}

Thus in general $H_{\varphi }^{\ast }$ takes values in the extended real
line. The problem is that the function $H_{\varphi }^{\ast }$ is in general
less well behaved than the function $h$ -- for example, if $h$ is
measurable/Baire, $H_{\varphi }^{\ast }$ need not be. The problem we address
here is the extent of this degradation -- saying \textit{exactly how much
less regular} than $h$ the limsup $H_{\varphi }^{\ast }$ may be. The nub is
the set $S$ on which $H_{\varphi }^{\ast }$ is finite. This set $S$ is an
additive semi-group on which the function $H_{\varphi }^{\ast }$ is
subadditive (see [BinO7]) -- or additive, if limits exist (see[BinO6]).
Furthermore, if $H$ has Borel graph then $H_{\varphi }^{\ast }\ $has $%
\mathbf{\Delta }_{2}^{1}$ graph (see below). But in the presence of certain
axioms of set-theory (for which see below) the $\mathbf{\Delta }_{2}^{1}$
sets have the Baire property and are measurable. Alternatively, if the $%
\mathbf{\Delta }_{2}^{1}$ character is witnessed by two $\mathbf{\Sigma }%
_{2}^{1}$ formulas $\Phi ,\Psi $ such that the equivalence%
\[
\Phi (x)\Longleftrightarrow \lnot \Psi (x)
\]%
is provable in $ZF,$ i.e. \textit{without reference} to $AC,$ then $A$ is
said to be \textit{provably} $\mathbf{\Delta }_{2}^{1}.$ It then turns out
that such sets are Baire/measurable -- see [FenN]. So in such circumstances
if $S$ is large in either of these two senses, then in fact $S$ contains a
half-line.

The extent of the degradation in passing from $h$ to $H_{\varphi }^{\ast }$
is addressed in the following result, which we call the First Character
Theorem, and then contrast it with two alternatives. These extend
corresponding results established in the Karamata context in [BinO8] and
differ from the former merely by duplicating assumptions previously made
only on $h$ there to identical ones on $\varphi .$

\bigskip

\noindent \textbf{Theorem 11 (First Character Theorem).} (i) \textit{If }$h$%
\textit{\ and }$\varphi $\textit{\ are Borel (have Borel graph), then the
graph of the function }%
\[
H^{\ast }(x)=\lim \sup\nolimits_{t\rightarrow \infty }[h(t+x\varphi
(t))-h(t)]
\]%
\textit{\ is a difference of two analytic sets, hence is measurable and }$%
\mathbf{\Delta }_{2}^{1}$. \textit{If the graphs of }$h$\textit{\ and }$%
\varphi $\textit{\ are }$\mathcal{F}_{\sigma }$\textit{, then the graph of }$%
H^{\ast }(x)$\textit{\ is Borel.}

\noindent (ii) \textit{If }$h$\textit{\ and }$\varphi $\textit{\ are
analytic (have analytic graph), then the graph of the function }$H^{\ast
}(x) $\textit{\ is }$\mathbf{\Pi }_{2}^{1}$\textit{.}

\noindent (iii) \textit{If }$h$\textit{\ and }$\varphi $\textit{\ are
co-analytic (have co-analytic graph), then the graph of the function }$%
H^{\ast }(x)$\textit{\ is }$\mathbf{\Pi }_{3}^{1}$\textit{.}

\bigskip

The next two results assume much more, in requiring the existence of a limit
(Th. 12) or a limit modulo an ultrafilter (Th. 13).

\bigskip

\noindent \textbf{Theorem 12 (Second Character Theorem).} \textit{If the
following limit exists:}%
\[
K_{h}(x):=\lim_{t\rightarrow \infty }[h(t+x\varphi (t))-h(t)],
\]%
\textit{and }$h,\varphi \in \mathbf{\Delta }_{2}^{1}$\ \textit{-- then the
graph of }$K_{h}$\textit{\ is }$\mathbf{\Delta }_{2}^{1}.$

\bigskip

\noindent \textbf{Theorem 13 (Third Character Theorem).} \textit{If the
function }$h$\textit{\ and the ultrafilter }$\mathcal{U}$\textit{\ (both on }%
$\omega $) \textit{are of class }$\mathbf{\Delta }_{2}^{1}$ -- \textit{then
so is: }%
\[
K_{h}^{\mathcal{U}}h(t):=\mathcal{U}\text{-}\lim_{n}[h(n+t\varphi
(n))-h(n)].
\]

The proofs of all three character theorems closely follow the proofs of the
Karamata special case in [BinO8, \S 4], by using just two amendment
procedures. Firstly, apply a \textit{replacement rule:} all uses of the
formula $y=h(x,t):=h(x+t)-h(t)$ ($h$ as there) are to be replaced by a
formalized conjunction of $y=h(x,s,t):=h(x+ts)-h(t)$ and $s=\varphi (x),$ as
follows. Translate these two formulas to `$(x,s,t,y)\in h$ \& $(x,s)\in
\varphi $' (interpreting $h$ and $\varphi $ as naming the graphs of the two
functions), and replace each $(x,t,y)\in h$ there by the the translate just
indicated here above. Secondly, apply an \textit{insertion rule}: insert the
variable $s$ everywhere to precede the variable $w$. An example of the
translation will suffice; here is a sample amendment:
\[
y=h(t+xs)-h(t)\Leftrightarrow (\exists u,v,s,w\in \mathbb{R}%
)r(x,t,y,u,v,s,w),
\]%
where $r(x,t,y,u,v,s,w)$ stands for:%
\begin{equation}
\lbrack \text{ }y=u-v\text{ \& }w=t+xs\text{ \& }(w,u)\in h\text{ \& }%
(t,v)\in h\text{ \& }(x,s)\in \varphi ].  \label{r}
\end{equation}

\noindent \textbf{Comment 1. }In the first theorem (as also in [BinO8]) we
deal with $H^{\ast }(x)=K_{h}^{\ast }(x):=\lim \sup_{t\rightarrow \infty
}\Delta _{x}^{\varphi }h(t).$ The results are also true for $\lim \sup
\Delta _{x}^{\varphi }h/\varphi (t)$ or $\lim \sup \Delta _{x}^{\varphi
}h/\psi (t).$ The proofs are essentially the same; one needs the same
assumptions on $\varphi $ (or $\psi $) as on $h.$

\noindent \textbf{Comment 2. }The last of the three theorems applies under
the assumption of G\"{o}del's Axiom $V=L$ (see [Dev, \textbf{\S }B.5,
453-489]), under which $\mathbf{\Delta _{2}^{1}}$ ultrafilters exist on $%
\omega $ (e.g. for Ramsey ultrafilters -- see [Z]). Above sets of natural
numbers are identified with real numbers (via indicator functions), and so
ultrafilters are subsets of $\mathbb{R}$ -- for background see [CoN], or
[HinS]. Th. 12 offers a midway position between the First and Second
Character Theorems.

In Th. 13 $K_{h}^{\mathcal{U}}h(t)$ is additive, whereas in Th. 11 one has
only sub-additivity (cf. BGT\ p. 62 equation (2.0.3)).

\noindent \textbf{Comment 3. }Replacing $h(n+t\varphi (n))-h(n)$ by $%
h(x(n)+t\varphi (x(n))-h(x(n))$, as in the Equivalence Theorem of [BinO3],
to take limits along a specified sequence\textit{\ }$\mathbf{x}:\omega
\rightarrow \omega ^{\omega },$ gives an `effective' version of the
character theorems -- given an effective descriptive character of $\mathbf{x.%
}$

\bigskip

\noindent \textbf{Acknowledgements. }We are most grateful to the Referee for
his extremely detailed, scholarly and helpful report, which has led to many
improvements. We also thank Guus Balkema for helpful comments.

\bigskip

\begin{center}
\textbf{References}
\end{center}

\noindent \lbrack Bec] A. Beck, \textsl{Continuous flows on the plane},
Grundl. math. Wiss. \textbf{201}, Springer, 1974.\newline
\noindent \lbrack Bin] N. H. Bingham, Riesz means and Beurling moving
averages, \textsl{Risk \& Stochatics} (Ragnar Norberg Festschrift, ed. P. M.
Barrieu), Imp. Coll. Press, to appear; arXiv 1502.07494.\newline
\noindent \lbrack BinG1] N. H. Bingham, C. M. Goldie, Extensions of regular
variation. II. Representations and indices. \textsl{Proc. London Math. Soc.}
(3) \textbf{44} (1982), 497--534.\newline
\noindent \lbrack BinG2] N. H. Bingham, C. M. Goldie, On one-sided Tauberian
conditions. \textsl{Analysis} \textbf{3} (1983), 159--188.\newline
\noindent \lbrack BinG3] N. H. Bingham, C. M. Goldie, Riesz means and
self-neglecting functions. \textsl{Math. Z.,} \textbf{199} (1988), 443--454.%
\newline
\noindent \lbrack BinGT] N. H. Bingham, C. M. Goldie and J. L. Teugels,
\textsl{Regular variation}, 2nd ed., Cambridge University Press, 1989 (1st
ed. 1987).\newline
\noindent \lbrack BinO1] N. H. Bingham and A. J. Ostaszewski, Beyond
Lebesgue and Baire: generic regular variation. \textsl{Coll. Math.} \textbf{%
116} (2009), 119-138.\newline
\noindent \lbrack BinO2] N. H. Bingham and A. J. Ostaszewski, The index
theorem of topological regular variation and its applications. \textsl{J.
Math. Anal. Appl.} \textbf{358} (2009), 238-248. \newline
\noindent \lbrack BinO3] N. H. Bingham and A. J. Ostaszewski, Infinite
combinatorics and the foundations of regular variation, \textsl{J. Math.
Anal. Appl.} \textbf{360} (2009), 518-529. \newline
\noindent \lbrack BinO4] N. H. Bingham and A. J. Ostaszewski, Beyond
Lebesgue and Baire II: Bitopology and measure-category duality. \textsl{%
Coll. Math.}, \textbf{121} (2010), 225-238.\newline
\noindent \lbrack BinO5] N. H. Bingham and A. J. Ostaszewski, Topological
regular variation. I: Slow variation; II: The fundamental theorems; III:
Regular variation. \textsl{Topology Appl.} \textbf{157} (2010), 1999-2013;
2014-2023; 2024-2037. \newline
\noindent \lbrack BinO6] N. H. Bingham and A. J. Ostaszewski, Normed groups:
Dichotomy and duality. \textsl{Dissertationes Math.} \textbf{472} (2010),
138p. \newline
\noindent \lbrack BinO7] N. H. Bingham and A. J. Ostaszewski, Kingman,
category and combinatorics. \textsl{Probability and Mathematical Genetics}
(Sir John Kingman Festschrift, ed. N. H. Bingham and C. M. Goldie), 135-168,
London Math. Soc. Lecture Notes in Mathematics \textbf{378}, CUP, 2010.%
\newline
\noindent \lbrack BinO8] N. H. Bingham and A. J. Ostaszewski: Regular
variation without limits, \textsl{J. Math. Anal. Appl.}, \textbf{370}
(2010), 322-338.\newline
\noindent \lbrack BinO9] N. H. Bingham and A. J. Ostaszewski, Dichotomy and
infinite combinatorics: the theorems of Steinhaus and Ostrowski. \textsl{%
Math. Proc. Camb. Phil. Soc.} \textbf{150} (2011), 1-22. \newline
\noindent \lbrack BinO10] N. H. Bingham and A. J. Ostaszewski: Beurling slow
and regular variation, \textsl{Trans. London. Math. Soc., }\textbf{1}
(2014), 29-56 (see also Part I: arXiv:1301.5894, Part II: arXiv:1307.5305).%
\newline
\noindent \lbrack BinO11] N. H. Bingham and A. J. Ostaszewski, Cauchy's
functional equation and extensions: Goldie's equation and inequality, the Go%
\l \k{a}b-Schinzel equation and Beurling's equation, \textsl{Aequationes
Math.}, \textbf{89} (2015), 1293-1310, arXiv1405.3947.\newline
\noindent \lbrack BinO12] N. H. Bingham and A. J. Ostaszewski: Additivity,
subadditivity and linearity: automatic continuity and quantifier weakening,
arXiv.1405.3948.\newline
\noindent \lbrack BinO13] N. H. Bingham and A. J. Ostaszewski:
Category-measure duality, Jensen convexity and Berz sublinearity, in
preparation.\textbf{\newline
}\noindent \lbrack Blo] S. Bloom, A characterization of B-slowly varying
functions. \textsl{Proc. Amer. Math. Soc.} \textbf{54} (1976), 243-250.
\newline
\noindent \lbrack Boa] R. P. Boas, \textsl{A primer of real functions}. 3rd
ed. Carus Math. Monographs 13, Math. Assoc. America, 1981. \newline
\noindent \lbrack BojK] R. Bojani\'{c} and J. Karamata, \textsl{On a class
of functions of regular asymptotic behavior, }Math. Research Center Tech.
Report 436, Madison, Wis. 1963; reprinted in \textsl{Selected papers of
Jovan Karamata} (ed. V. Mari\'{c}, Zevod za Ud\v{z}benike, Beograd, 2009),
545-569.\newline
\noindent \lbrack Brz1] J. Brzd\k{e}k, The Go\l \k{a}b-Schinzel equation and
its generalizations, \textsl{Aequat. Math.} \textbf{70} (2005), 14-24.%
\newline
\noindent \lbrack Brz2] J. Brzd\k{e}k, A remark on solutions of a
generalization of the addition formulae, \textsl{Aequationes Math.}, \textbf{%
71} (2006), 288--293.\newline
\noindent \lbrack BrzM] J. Brzd\k{e}k and A. Mure\'{n}ko, On a conditional Go%
\l \k{a}b-Schinzel equation, \textsl{Arch. Math.} \textbf{84} (2005),
503-511.\newline
\noindent \lbrack Chu] J. Chudziak, Semigroup-valued solutions of the Go\l
\k{a}b-Schinzel type functional equation, \textsl{Abh. Math. Sem. Univ.
Hamburg,} \textbf{76} (2006), 91-98.\newline
\noindent \lbrack CoN] W.W. Comfort, S. Negrepontis, \textsl{The theory of
ultrafilters.} Die Grundlehren der mathematischen Wissenschaften, Band
\textbf{211}. Springer-Verlag, New York-Heidelberg, 1974.\newline
\noindent \lbrack Dev] K. J. Devlin, \textsl{Constructibility}, Springer
1984.\newline
\noindent \lbrack Dra] F. R. Drake, D. Singh, \textsl{Intermediate set
theory.} Wiley, 1996\newline
\noindent \lbrack FenN] J. E. Fenstad, D. Normann, On absolutely measurable
sets. \textsl{Fund. Math.} \textbf{81} (1973/74), no. 2, 91--98.\newline
\noindent \lbrack dH] L. de Haan, On regular variation and its applications
to the weak convergence of sample extremes. Math. Centre Tracts \textbf{32},
Amsterdam 1970.\newline
\noindent \lbrack HinS] N. Hindman, D. Strauss, \textsl{Algebra in the Stone-%
\v{C}ech compactification. Theory and applications.} 2nd rev. ed., de
Gruyter, 2012. (1st. ed. 1998)\newline
\noindent \lbrack Hob] E.W. Hobson, \textsl{The theory of functions of a
real variable and the theory of Fourier's Series}, Vol. 2, 2$^{\text{nd}}$
ed., CUP, 1926.\newline
\noindent \lbrack Jav] P. Javor, On the general solution of the functional
equation $f(x+yf(x))=f(x)f(y).$ \textsl{Aequat. Math.} \textbf{1} (1968),
235-238.\newline
\noindent \lbrack JayR] J. Jayne and C. A. Rogers, \textsl{K-analytic sets, }%
Part 1 (p.1-181) in [Rog].\newline
\noindent \lbrack Kec] A. S. Kechris: \textsl{Classical Descriptive Set
Theory.} Grad. Texts in Math. \textbf{156}, Springer, 1995.\newline
\noindent \lbrack KliW] J. Klippert, G. Williams, Uniform convergence of a
sequence of functions at a point, \textsl{Int. J. Math. Ed. in Sc. and Tech.}%
, \textbf{33.1} (2002), 51-58.\newline
\noindent \lbrack Kor] J. Korevaar, \textsl{Tauberian theorems: A century of
development}. Grundl. math. Wiss. \textbf{329}, Springer, 2004.\newline
\noindent \lbrack Kuc] M. Kuczma, \textsl{An introduction to the theory of
functional equations and inequalities. Cauchy's equation and Jensen's
inequality.} 2nd ed., Birkh\"{a}user, 2009 [1st ed. PWN, Warszawa, 1985].%
\newline
\noindent \lbrack Kun] K. Kunen, \textsl{Set theory. An introduction to
independence proofs.} Reprint of the 1980 original. Studies in Logic and the
Foundations of Mathematics \textbf{102}. North-Holland, 1983.\newline
\noindent \lbrack MaKe] D. A. Martin, A.S. Kechris, Infinite games and
effective descriptive set theory, Part 4 (p. 403-470) in [Rog].\newline
\noindent \lbrack Mos] Y. N. Moschovakis{, \textsl{Descriptive set theory,}
Studies in Logic and the Foundations of Math. }\textbf{100},{\
North-Holland, Amsterdam, 1980.}\newline
\noindent \lbrack Mur] A. Mure\'{n}ko, On the general solution of a
generalization of the Go\l \k{a}b-Schinzel equation, \textsl{Aequat. Math.},
\textbf{77} (2009), 107-118.\newline
\noindent \lbrack MySw] J. Mycielski and S. \'{S}wierczkowski, On the
Lebesgue measurability and the axiom of determinateness, \textsl{Fund. Math.}
\textbf{54} (1964), 67--71.\newline
\noindent \lbrack Ost1] A. J. Ostaszewski, Regular variation, topological
dynamics, and the Uniform Boundedness Theorem, \textsl{Top. Proc.}, \textbf{%
36} (2010), 305-336.\newline
\noindent \lbrack Ost2] A. J. Ostaszewski, Beyond Lebesgue and Baire III:
Steinhaus' Theorem and its descendants, \textsl{Topology Appl.} \textbf{160}
(2013), 1144-1154.\newline
\noindent \lbrack Ost3] A.J. Ostaszewski, Beurling regular variation, Bloom
dichotomy, and the Go\l \k{a}b-Schinzel functional equation, \textsl{%
Aequationes Math.} \textbf{89} (2015), 725-744.\newline
\noindent \lbrack Ost4] A. J. Ostaszewski, Homomorphisms from functional
equations: The Goldie Equation, \textsl{Aequationes Math.}, DOI
10.1007/s00010-015-0357-z, arXiv.org/abs/1407.4089.\newline
\noindent \lbrack Oxt] J. C. Oxtoby: \textsl{Measure and category}, 2nd ed.
Graduate Texts in Math. \textbf{2}, Springer, 1980.\newline
\noindent \lbrack PolS] G. P\'{o}lya and G. Szeg\H{o}, \textsl{Aufgaben und
Lehrs\"{a}tze aus der Analysis} Vol. I. Grundl. math. Wiss. XIX, Springer,
1925. \newline
\noindent \lbrack Pop] C. G. Popa, Sur l'\'{e}quation fonctionelle $%
f[x+yf(x)]=f(x)f(y),$ \textsl{Ann. Polon. Math.} \textbf{17} (1965), 193-198.%
\newline
\noindent \lbrack Rog] C. A. Rogers, J. Jayne, C. Dellacherie, F. Tops\o e,
J. Hoffmann-J\o rgensen, D. A. Martin, A. S. Kechris, A. H. Stone, \textsl{%
Analytic sets,} Academic Press,1980.\newline
\noindent \lbrack Rud] W. Rudin, \textsl{Principles of Mathematical Analysis}%
. 3rd ed. McGraw-Hill, 1976 (1st ed. 1953).\newline
\noindent \lbrack Ste] H. Stetkaer, \textsl{Functional equations on groups},
World Scientific, 2013.\newline
\noindent \lbrack Wid] D. V. Widder, \textsl{The Laplace Transform},
Princeton, 1972.\newline
\noindent \lbrack Wie] N. Wiener, \textsl{The Fourier integral \& certain of
its applications}, CUP, 1988.

\noindent \lbrack Z] J. Zapletal, Terminal Notions, \textsl{Bull. Symbolic
Logic} \textbf{5} (1999), 470-478.\newline

\bigskip

\noindent Mathematics Department, Imperial College, London SW7 2AZ;
n.bingham@ic.ac.uk \newline
Mathematics Department, London School of Economics, Houghton Street, London
WC2A 2AE; A.J.Ostaszewski@lse.ac.uk

\bigskip

\section*{Appendix: Global bounds}

Below we need Bloom's [Blo] result that for $x$ large enough the Beck
sequence $x_{n}^{u}$ defined recursively by its starting value $x$ and the
step-size $u:$%
\[
x_{n+1}^{u}=x_{n}\circ _{\varphi }u=x_{n}+u\varphi (x_{n}),\text{ with }%
x_{0}^{u}=x,\text{ }x_{1}^{u}=x+u\varphi (x)
\]%
is divergent (see [BinO-B, \S 9] and compare [Ost-B, \S 6]). Say for $x\geq
x_{0}.$

We briefly review a number of examples of Beck sequences; Example 2 is
crucial.

\bigskip

\noindent \textbf{Example }1. $a_{\varphi }^{n}=a_{n}^{a},$ (with $a_{n}^{a}$
as above) so that $a_{\varphi }^{n+1}=a_{\varphi }^{n}\circ _{\varphi
}a=a_{\varphi }^{n}+a\varphi (a_{\varphi }^{n}).$ Performing the recurrence
the other way about, $u_{n+1}=u\circ _{\varphi }u_{n}=u+u_{n}\varphi (u)$
generates a GP:%
\[
u_{n}=(1-\varphi (u)^{n+1})\cdot u/(1-\varphi (u)),
\]%
with%
\[
u_{n+1}-u_{n}=(u_{n}-u_{n-1})\varphi (u)=...=u\varphi (u)^{n}.
\]%
For $\varphi \in GS$ the two are the same. They are not altogether
dissimilar, as the other one has%
\[
a_{\varphi }^{k}=a[1+\varphi (a)+\varphi (a_{\varphi }^{2})+...+\varphi
(a_{\varphi }^{k-1})],
\]%
and, assuming divergence, the term-on-term growth is%
\[
\varphi (a_{\varphi }^{k})/\varphi (a_{\varphi }^{k-1})=\varphi (a_{\varphi
}^{k-1}+a\varphi (a_{\varphi }^{k-1}))/\varphi (a_{\varphi
}^{k-1})\rightarrow \eta ^{\varphi }(a),
\]%
so the \textit{series} behaves, up to a multiplier $\varphi (a_{\varphi
}^{k}),$ eventually like
\[
\sum_{j<k}\eta ^{\varphi }(a)^{j}=(1-\eta ^{\varphi }(a_{\eta
}^{k}))/(1-\eta ^{\varphi }(a)).
\]%
\noindent \textbf{Example 2. }Consider the sequence%
\[
a_{\varphi x}^{n+1}:=a_{\varphi x}^{n}\circ _{\varphi x}a=a_{\varphi
x}^{n}+a\eta _{x}^{\varphi }(a_{\varphi x}^{n})\text{ with }a_{\varphi
x}^{1}=a,
\]%
where $a$ is fixed; on the back of Example 1 we guess that since uniformly
in $x$
\[
\eta _{x}^{\varphi }(a)\rightarrow \eta _{\rho }(a),
\]%
this $a_{\varphi x}^{n}$ is a divergent sequence for $x$ large enough, say $%
x>x_{a}$. Indeed, it is -- see the proof of Prop. 11; this is to be expected
from the related iteration%
\[
a_{\eta }^{n+1}:=a_{\eta }^{n}\circ _{\eta }a=a_{\eta }^{n}+a\eta _{\rho
}(a_{\eta }^{n})\text{ with }a_{\eta }^{1}=a,
\]%
where for $\rho =0$ growth is linear: $\eta (a_{\eta }^{n})=na$, whereas for
$\rho >0$ it is exponential:%
\[
\eta (a_{\eta }^{n})=\eta (a_{\eta }^{n-1}\circ _{\eta }a)=\eta (a_{\eta
}^{n-1})\eta (a)=...=\eta _{\rho }(a)^{n}=(1+\rho a)^{n}.
\]

Below we need the solution of a recurrence; we present this as a lemma,
delaying the calculation to the end.

\bigskip

\noindent \textbf{Lemma 3}. \textit{The solution of }$bv_{n+1}-v_{n}=r^{n}$%
\textit{\ for }$br\neq 1$\textit{\ is}%
\begin{equation}
v_{n}=r^{n}/(br-1)+b^{1-n}(v_{1}-r/(br-1)).  \tag{soln}
\end{equation}%
\textit{If }$b=\eta _{\rho }(a)$ \textit{with }$\rho >0,$\textit{\ }$%
v_{1}=1/(\rho a),$\textit{\ }$r=1\pm \delta ,$\textit{\ with }$\delta
=\varepsilon \rho a/\eta _{\rho }(a)$ \textit{and }$0<\varepsilon <1$\textit{%
, then}%
\[
v_{1}-r/(br-1)=\frac{\varepsilon /(\eta _{\rho }(a)\rho a)}{(1+\varepsilon )}%
\text{ or}-\frac{\varepsilon /(\eta _{\rho }(a)\rho a)}{(1-\varepsilon )}%
\text{.}
\]

We now proceed to verify the details of Prop. 11 in \S 10.

\bigskip

\noindent \textbf{Proof of Prop. 11.} Fix $a,\rho >0$ and $0<\varepsilon <1.$
Taking $\delta :=\varepsilon \rho a/\eta (a),$
\[
\eta (a)\pm \rho a\varepsilon =(1+\rho a(1\pm \varepsilon ))=\eta (a(1\pm
\varepsilon ))=\eta (a)(1\pm \delta ).
\]%
In particular, $\eta (a)(1-\delta )=\eta (a(1-\varepsilon ))>1,$ since $%
\varepsilon <1.$ Since $\eta _{x}(a)\rightarrow \eta (a),$ there is $%
X=X_{a,\varepsilon }$ with%
\begin{equation}
|\eta (a)-\eta _{x}(a)|<\rho a\varepsilon :\qquad \eta (a)(1-\delta )<\eta
_{x}(a)<\eta (a)(1+\delta )\qquad (x>X).  \tag{$\delta $-bd}
\end{equation}%
\noindent (i) By Prop. 2(v), for $y_{i}$ running through $x\circ _{\varphi
}a_{\varphi x}^{m-1},$ $x\circ _{\varphi }a_{\varphi x}^{m-2},...,x>X,$%
\begin{equation}
\eta _{x}(a_{\varphi x}^{m})=\prod\nolimits_{i=1}^{m}\eta _{y_{i}}(a),
\tag{prod}
\end{equation}%
so that, by ($\delta $-bd),
\[
\eta (a(1-\varepsilon ))\leq \eta _{x}(a_{\varphi x}^{m})^{1/m}\leq \eta
(a(1+\varepsilon )).
\]%
\noindent (ii) As $\eta \in GS,$ $\eta (a_{\varphi x}^{n+1})=\eta
(a_{\varphi x}^{n}+a\eta _{x}(a_{\varphi x}^{n}))=\eta (a_{\varphi
x}^{n})\eta (a\eta _{x}(a_{\varphi x}^{n})/\eta (a_{\varphi x}^{n})).$ So%
\[
\eta (a_{\varphi x}^{n+1})/\eta (a_{\varphi x}^{n})=1+\rho a\eta
_{x}(a_{\varphi x}^{n})/\eta (a_{\varphi x}^{n}):\qquad \eta (a_{\varphi
x}^{n+1})-\eta (a_{\varphi x}^{n})=\rho a\eta _{x}(a_{\varphi x}^{n}).
\]%
Putting $u_{n}:=\eta (a_{\varphi x}^{n})/(\rho a\eta (a)^{n}),$ so that $%
u_{1}=1/(\rho a),$ and using ($\delta $-bd) again,
\[
(1-\delta )^{n}\leq \frac{\eta (a_{\varphi x}^{n+1})-\eta (a_{\varphi x}^{n})%
}{\rho a\eta (a)^{n}}=\eta (a)u_{n+1}-u_{n}\leq (1+\delta )^{n}.
\]%
As $\eta (a)(1\pm \delta )\neq 1$, apply Lemma 3 to $b=\eta (a)$ and $r=1\pm
\delta ;$ then%
\[
\frac{(1-\delta )^{n}\eta (a)^{n}}{1-\varepsilon }-\frac{\varepsilon }{%
1-\varepsilon }\leq \eta (a_{\varphi x}^{n})\leq \frac{(1+\delta )^{n}\eta
(a)^{n}}{1+\varepsilon }+\frac{\varepsilon }{1+\varepsilon }.
\]%
\noindent (iii) As $\eta (a)(1-\delta )>1,$ the left inequality implies $%
a_{\varphi x}^{m}$ is \textit{divergent}.

\noindent (iv) If $a_{\varphi x}^{m}\leq u<a_{\varphi x}^{m+1},$ then (as $%
\eta _{\rho }$ is monotone), $\eta (a_{\varphi x}^{m})\leq 1+\rho u\leq \eta
(a_{\varphi x}^{m+1});$ so, for $x>X_{e}$%
\[
\frac{\eta (a(1-\varepsilon ))^{m}}{(1-\varepsilon )}-\frac{1}{1-\varepsilon
}\leq \rho u<\frac{\eta (a(1+\varepsilon ))^{m+1}}{1+\varepsilon }-\frac{1}{%
1+\varepsilon }.
\]%
So for $\varepsilon <1/2$
\[
\frac{\eta (a(1-\varepsilon ))^{m}}{(1-\varepsilon )}-2\leq \rho u<\frac{%
\eta (a(1+\varepsilon ))^{m+1}}{1+\varepsilon }.
\]%
So for $u>1,$ as $2+\rho u<(2+\rho )u,$%
\[
\frac{\eta (a(1-\varepsilon ))^{m}}{(1-\varepsilon )(2+\rho )}\leq u<\frac{%
\eta (a(1+\varepsilon ))^{m+1}}{\rho (1+\varepsilon )},
\]%
where $\log \eta (a(1-\varepsilon ))>0.$ Taking%
\[
C_{-}(\rho ,a,\varepsilon ):=\log \frac{\eta (a(1-\varepsilon ))}{(\rho
+2)(1-\varepsilon )},\text{ }C_{+}(\rho ,a,\varepsilon ):=\log \frac{\eta
(a(1+\varepsilon ))}{\rho (1+\varepsilon )},
\]%
\[
mC_{-}\leq \log u<(m+1)C_{+}\qquad (u\geq a>1\text{ \& }x\geq
X_{a,\varepsilon }).\qquad \square
\]

\bigskip

We are now ready to prove Th. 9 of \S 10.

\bigskip

\noindent \textbf{Proof of Theorem 10.} (This parallels the tail end of the
proof in BGT\ of Th. 2.0.1.) W.l.o.g. we assume that $\eta ^{\varphi
}(x)=1+\rho x$ with $\rho >0,$ as the case $\rho =0$ is already known. By
Theorem 8 (UBT) in \S 10, for any $a\geq a_{0}$%
\[
\limsup_{x\rightarrow \infty }\left( \sup_{a\leq u\leq 2a\eta (a)}h(x\circ
_{\varphi }u)-h(x)\right) <\infty .
\]%
So there is $C_{a}$ such that%
\[
\sup_{a\leq u\leq 2a\eta (a)}h(x\circ _{\varphi }u)-h(x)<C_{a},
\]%
for all large enough $x,$ say for $x>x_{a}.$ Choose $a>\max \{a_{0},x_{a}\}.$

As at the start of the proof of Prop. 11, but specializing to $\varepsilon
=1,$ take $\delta :=\rho a/\eta (a)$ to obtain ($\delta $-$bd$) for $x>X:$%
\begin{equation}
\eta (a)(1-\delta )<\eta _{x}(a)<\eta (a)(1+\delta )\qquad (x>X).
\tag{$\ast \ast $}
\end{equation}%
For $x>X,$ fix $u\geq a=a_{\varphi x}^{1}.$ Then, by Prop. 11(iii), we may
choose $m=m_{x}(u)$ such that
\[
a_{\varphi x}^{m-1}<a_{\varphi x}^{m}\leq u\leq a_{\varphi x}^{m+1}.
\]%
Now put $d:=(u-a_{\varphi x}^{m-1})/\eta _{x}(a_{\varphi x}^{m-1}),$ so that
$u=a_{\varphi x}^{m-1}\circ _{\varphi x}d;$ then
\[
x\circ _{\varphi }u=[x+a_{\varphi x}^{m-1}\varphi (x)]+d\varphi
(x+a_{\varphi x}^{m-1}\varphi (x))=y\circ _{\varphi }d,
\]%
with $y=x\circ _{\varphi }a_{\varphi x}^{m-1};$ referring to $[a_{\varphi
x}^{m}-a_{\varphi x}^{m-1}]+[a_{\varphi x}^{m+1}-a_{\varphi x}^{m}]=a\eta
_{x}(a_{\varphi x}^{m-1})+a\eta _{x}(a_{\varphi x}^{m})$ and to $%
u-a_{\varphi x}^{m-1}=d\eta _{x}(a_{\varphi x}^{m-1}),$
\[
a\eta _{x}(a_{\varphi x}^{m-1})\leq d\eta _{x}(a_{\varphi x}^{m-1})<a\eta
_{x}(a_{\varphi x}^{m-1})+a\eta _{x}(a_{\varphi x}^{m})=a\eta
_{x}(a_{\varphi x}^{m-1})+a\eta _{y}(a)\eta _{x}(a_{\varphi x}^{m-1}),
\]%
as in Prop. 2(v). But by (**) above, since $y\geq x>x_{a},$
\[
a\leq d<a(1+\eta _{y}(a))<a(1+\eta (a)(1+\delta ))<a(1+\rho a+\eta
(a))=2a\eta (a),
\]%
as $\delta \eta (a)=\rho a.$ So by choice of $C_{a},$
\[
h(x\circ _{\varphi }u)-h(x\circ _{\varphi }a_{\varphi x}^{m-1})=h(y\circ
_{\varphi }d)-h(y)<C_{a},
\]%
as $d\in \lbrack a,2a\eta (a)].$ As in Prop. 11,%
\[
x\circ _{\varphi }a_{\varphi x}^{n+1}=x\circ _{\varphi }(a_{\varphi
x}^{n}\circ _{\varphi }a)=(x\circ _{\varphi }a_{\varphi x}^{n})\circ
_{\varphi }a,
\]%
and, setting $y_{k}=x\circ _{\varphi }a_{\varphi x}^{k}$ for $k=0,...,m-1,$%
\[
h(x\circ _{\varphi }a_{\varphi x}^{k+1})-h(x\circ _{\varphi }a_{\varphi
x}^{k})=h((x\circ _{\varphi }a_{\varphi x}^{k})\circ _{\varphi }a)-h(x\circ
_{\varphi }a_{\varphi x}^{k})=h(y_{k}\circ _{\varphi }a)-h(y_{k})<C_{a},
\]%
since $y_{k}\geq x>x_{a}.$ So for $x>x_{a}$%
\[
h(x\circ _{\varphi }u)-h(x)=h(x\circ _{\varphi }u)-h(x\circ _{\varphi
}a_{\varphi x}^{m-1})+\sum\limits_{k=1}^{m-1}(h(x\circ _{\varphi }a_{\varphi
x}^{k})-h(x\circ _{\varphi }a_{\varphi x}^{k-1}))<mC_{a}.
\]%
Again by Prop. 11, there is a constant $C$ such that%
\[
m\leq C\log u.
\]%
Taking $K=C_{a}C$ yields the desired inequality. $\square $

\bigskip

\noindent \textbf{Proof of Lemma 3. }A particular solution is $r^{n}/(br-1),$
$bw_{n+1}-w_{n}=0$ for $w_{n}=v_{n}-r^{n}/(br-1)$ and $w_{n}=w_{1}b^{1-n},$
where $w_{1}=v_{1}-r/(br-1).$

For $b=\eta _{\rho }(a),$ $v_{1}=1/(\rho a)$ and $r=1\pm \delta ,$ we
calculate that%
\begin{eqnarray*}
\rho aw_{1} &=&\frac{[\eta (a)(1\pm \delta )-1]-\rho a(1\pm \delta )}{%
(1+\rho a)(1\pm \delta )-1}=\frac{[(1+\rho a)(1\pm \delta )-1]-\rho a(1\pm
\delta )}{\rho a+\eta (a)(\pm \delta )} \\
&=&\pm \frac{\delta }{\rho a+\eta (a)(\pm \delta )}=\pm \frac{\varepsilon
/\eta (a)}{(1+(\pm 1)\varepsilon )}=\frac{\varepsilon /\eta (a)}{%
(1+\varepsilon )},\text{ or }-\frac{\varepsilon /\eta (a)}{(1-\varepsilon )}%
\text{ (}-\text{) }.
\end{eqnarray*}

We close with

\bigskip

\noindent \textbf{Proof of the Characterization Theorem (Uniform
semicontinuity).. }In the notation above, for $n>m$%
\[
f(t)-\varepsilon \leq \inf \{f_{n}(s):s\in I_{\delta }(t)\}\leq \sup
\{f_{n}(s):s\in I_{\delta }(t)\}\leq f(t)+\varepsilon .
\]%
So%
\[
f(t)-\varepsilon \leq \lim \inf_{n}\inf \{f_{n}(s):s\in I_{\delta }(t)\}\leq
\lim \sup_{n}\sup \{f_{n}(s):s\in I_{\delta }(t)\}\leq f(t)+\varepsilon .
\]%
We my now take limits as $\delta \downarrow 0$ to obtain%
\[
f(t)-\varepsilon \leq \lim_{\delta \downarrow 0}\lim \inf_{n}\inf
\{f_{n}(s):s\in I_{\delta }(t)\}\leq \lim_{\delta \downarrow 0}\lim
\sup_{n}\sup \{f_{n}(s):s\in I_{\delta }(t)\}\leq f(t)+\varepsilon .
\]%
As $\varepsilon >0$ was arbitrary,%
\begin{eqnarray*}
f(t) &=&\lim_{\delta \downarrow 0}\lim \sup_{n}\sup \{f_{n}(s):s\in
I_{\delta }(t)\} \\
&=&\lim_{\delta \downarrow 0}\lim \inf_{n}\inf \{f_{n}(s):s\in I_{\delta
}(t)\}.
\end{eqnarray*}%
Now suppose that $f(t)=\lim_{\delta \downarrow 0}\lim \sup_{n}\sup
\{f_{n}(s):s\in I_{\delta }(t)\}$ and $\varepsilon >0.$ Then for some $%
\delta >0$%
\[
\lim \sup_{n}\sup \{f_{n}(s):s\in I_{\delta }(t)\}<f(t)+\varepsilon ,
\]%
and so there is $N_{t}$ such that for $n>N_{t}$%
\[
\sup \{f_{n}(s):s\in I_{\delta }(t)\}<f(t)+\varepsilon
\]%
and so
\[
f_{n}(s)<f(t)+\varepsilon \text{ for }n>N_{t}\text{ and }s\in I_{\delta
}(t).
\]%
By a similar argument there is $\delta ^{\prime }$ and $N_{t}^{\prime }$ so
that%
\[
f_{n}(s)>f(t)-\varepsilon \text{ for }n>N_{t}^{\prime }\text{ and }s\in
I_{\delta ^{\prime }}(t).\text{ }\square
\]

\end{document}